\documentclass[tre,nonblindrev]{informs3}

\OneAndAHalfSpacedXI

\usepackage{natbib}
 \bibpunct[, ]{(}{)}{,}{a}{}{,}
 \def\bibfont{\small}

\makeatletter

\newcommand*{\rom}[1]{\expandafter\@slowromancap\romannumeral #1@}
\makeatother
\usepackage[utf8]{inputenc}
\usepackage{float}
\usepackage{subcaption}
\usepackage{tabularx}
\usepackage{adjustbox}
\usepackage{lineno}

\usepackage{amsmath}
\usepackage{amssymb}
\usepackage{breqn}
\usepackage[linesnumbered,ruled,vlined]{algorithm2e}
\usepackage[flushleft]{threeparttable}
\usepackage{bookmark}
\usepackage{tikz}
\usepackage{xcolor}
\usepackage{eurosym}
\usepackage{longtable}
\usepackage{xltabular}
\usepackage{booktabs}
\usepackage{hyperref}

\usepackage{appendix}
\usepackage{pdflscape}
\usepackage{afterpage}
\usepackage{capt-of}
\usepackage{lipsum}
\usepackage{rotating}
\usepackage{doi}
\usepackage{rotating}
\usepackage{multirow}

\tikzstyle{startstop} = [rectangle,rounded corners, minimum width=3cm,minimum height=1cm,text centered, draw=black,fill=red!30]
\tikzstyle{io} = [trapezium, trapezium left angle = 70,trapezium right angle=110,minimum width=3cm,minimum height=1cm,text centered,draw=black,fill=blue!30]
\tikzstyle{process} = [rectangle,minimum width=3cm,minimum height=1cm,text centered,text width =3cm,draw=black,fill=orange!30]
\tikzstyle{decision} = [diamond,minimum width=3cm,minimum height=1cm,text centered,draw=black,fill=green!30]
\tikzstyle{arrow} = [thick,->,>=stealth]
\tikzstyle{point}=[coordinate,on grid,]
\tikzstyle{block} = [draw, fill=blue!10, chamfered rectangle, minimum height=3em, minimum width=6em]

\usepackage{graphicx}
\usepackage{longtable,tabu}

\usetikzlibrary{shapes,arrows,chains}

\setcitestyle{authoryear,open={(},close={)}}

\TheoremsNumberedThrough

\EquationsNumberedThrough

\usepackage{natbib}
 \bibpunct[, ]{(}{)}{,}{a}{}{,}
 \def\bibfont{\small}

\begin{document}


\RUNTITLE{Dynamic Preference-based Multi-modal Trip Planning of Public Transport and Shared Mobility}

\RUNAUTHOR{Zhang, Oded, and Azadeh.}

\TITLE{{Dynamic Preference-based Multi-modal Trip Planning of Public Transport and Shared Mobility}}

\ARTICLEAUTHORS{
\AUTHOR{Yimeng Zhang$^{1,2}$\footnote{Corresponding author}, Oded Cats$^1$, Shadi Sharif Azadeh$^1$}
\AFF{\textsuperscript{1}Department of Transport \& Planning, Delft University of Technology, Delft, The Netherlands, \\
\textsuperscript{2}School of Transportation \& Logistics, Southwest Jiaotong University, Chengdu, China\\
\EMAIL{Yimeng.Zhang@tudelft.nl}}

}

\ABSTRACT{The shift from private vehicles to public and shared transport is crucial to reducing emissions and meeting climate targets. Consequently, there is an urgent need to develop a multimodal transport trip planning approach that integrates public transport and shared mobility solutions, offering viable alternatives to private vehicle use. To this end, we propose a preference-based optimization framework for multi-modal trip planning with public transport, ride-pooling services, and shared micro-mobility fleets. We introduce a mixed-integer programming model that incorporates preferences into the objective function of the mathematical model. We present a meta-heuristic framework that incorporates a customized Adaptive Large Neighborhood Search algorithm and other tailored algorithms, to effectively manage dynamic requests through a rolling horizon approach. Numerical experiments are conducted using real transport network data in a suburban area of Rotterdam, the Netherlands. Model application results demonstrate that the proposed algorithm can efficiently obtain near-optimal solutions. Managerial insights are gained from comprehensive experiments that consider various passenger segments, costs of micro-mobility vehicles, and availability fluctuation of shared mobility.}


\maketitle

\section{Introduction}\label{intro}

{By 2050, urban populations are projected to rise to 82\%, intensifying the need for sustainable urban mobility solutions \citep{World_Urbanization_Prospects}. Currently, private vehicles dominate road space and energy consumption, using over five times the space \citep{rode2017accessibility} and six times the energy \citep{Does_Bus_Transit_Reduce} compared to public transport or shared modes, contributing significantly to urban congestion, pollution, and greenhouse gas emissions. The shift from private vehicles to public transport and shared mobility is essential to meet climate targets, particularly as road transport accounts for nearly 70\% of transport-related emissions in Europe \citep{A_European_Strategy_low_emission_mobility}, 79\% in China \citep{pei2023analysis}, and 72.8\% in the United States \citep{US_emissions}.}

{In urban transport, except for private vehicles, a diverse array of transport services can potentially cater to the diverse needs of heterogeneous travelers. Passengers select either a single service or a combination thereof in order to fulfill their travel needs depending on their preferences and service availability. To foster better match between passengers and services, innovative platforms have emerged, such as Mobility as a Service (MaaS) \citep{wong2020mobility}. The interaction between passengers and such a platform is illustrated in Figure \ref{fig:relationships among passengers, service provider, and services}.
One essential task of the platform is to return multi-modal trip plans for passengers in response to travel requests while considering the optimal allocation of services based on the current demand patterns.} While there have been pioneering investigations into the multi-modal transport planning of Public Transport (PT) and ride-sharing services \citep{yu2021optimal,molenbruch2021analyzing}, the integration of micro-mobility solutions into PT has only recently begun to receive significant attention \citep{zhu2020understanding,liang2024dynamic}.
In recent years, there has been a notable increase in the number of policy initiatives specifically designed to promote the widespread adoption of micro-mobility services as a means of achieving environmentally sustainable transport \citep{abduljabbar2021role}. Furthermore, existing research predominantly concentrates on optimizing the supply aspect, {such as finding the optimal assignment of vehicles} \citep{ma2017demand}, or delves into demand analysis, typically employing choice models \citep{vij2020consumer, ho2020public}. {Addressing preference-based multi-modal trip planning of PT and shared mobility remains less explored, especially for heterogeneous users \citep{azadeh2022choice}.} The lack of consideration for individual preferences means that services may be designed as ``one-size-fits-all", potentially leading to reduced service quality. {Without incorporating preferences, trip planning may not prioritize the individual needs of users, {which affects the service level} \citep{zhang2022synchromodalb}, thereby resulting in unattractive travel offers. This could result in inconvenient travel experiences, potentially discouraging users from adopting multi-modal transport.}

\begin{figure}
    \centering
    \includegraphics[scale=0.6]{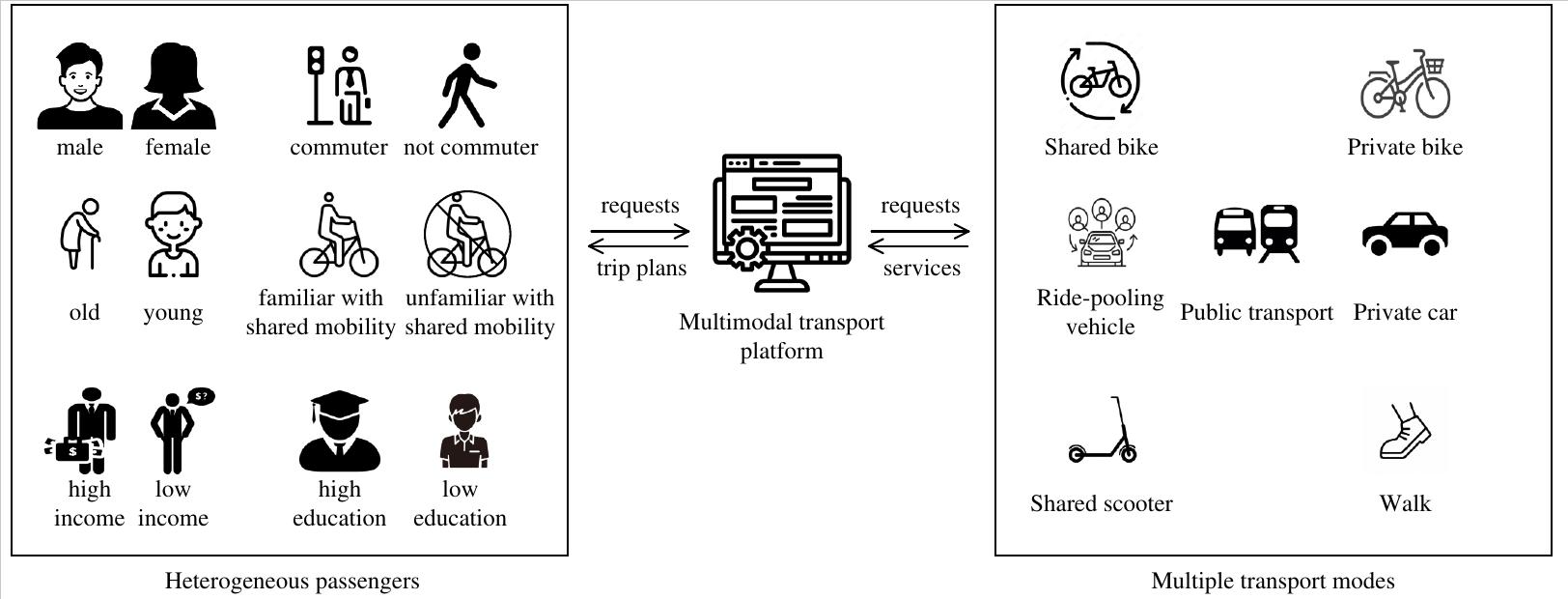}
    \caption{Interaction between passengers and a mobility platform}
    \label{fig:relationships among passengers, service provider, and services}
\end{figure}

{The preference-based multi-modal trip planning remains relatively unexplored, primarily due to its inherent complexities:
(a) Modal Integration: {Coordinating multiple transport modes seamlessly is a challenging task, which involves synchronization of schedules with integrated planning approaches.} (b) Preference Integration: Integrating the preferences of passengers into the planning process amplifies the planning complexity. Passengers have diverse and often conflicting preferences. Some might prioritize speed and direct routes, while others may prefer more economical or environmentally friendly options. This heterogeneity requires the platform to consider a wide range of preferences, making the planning process more intricate. To incorporate individual preferences, the platform needs to customize services and routes, leading to a multitude of potential combinations. This customization can result in a vast number of possibilities that must be evaluated, making the planning process computationally intensive \citep{azadeh2022choiceb}. To handle these complexities, the platform may need to develop routing and scheduling algorithms that can efficiently match passengers with their preferred services.}

 {To reduce private vehicle usage and its associated emissions, we propose a tailored preference-based Mixed Integer Programming (MIP) model and a meta-heuristic approach for multi-modal trip planning that integrates public transport and shared mobility, offering passengers a competitive and sustainable alternative.} Our focus extends to encompass line- and schedule-based PT modes, such as the metro, train, and bus, alongside Shared Mobility (SM) options, including ride-pooling and micro-mobility services. Our contributions include: (a) {{Introducing} a preference-based multi-modal trip planning framework, which seamlessly integrates PT and SM services. (b) {{Proposing} a dynamic planning approach based on a rolling horizon framework, which provides multimodal trip plans to dynamic passenger requests.}} (c) {Providing} valuable managerial insights through extensive numerical experiments conducted on a real-world transport network and different customer segments. These experiments include scenarios involving alterations in the availability and costs of SM, providing comprehensive insights into the adoption of multi-modal transport for a case study in the urban agglomeration area of Rotterdam, the Netherlands.

{The rest of this paper is structured as follows: Section \ref{Literature Review} offers a concise literature review. Section \ref{Problem description} describes the preference-based multi-modal trip planning problem addressed in this study. In Section \ref{Mathematical model for preference-based multimodal trip planning}, the mathematical model formulation is presented. Section \ref{Solution approach for preference-based multimodal trip planning} details our solution approach, including the rolling horizon approach, the integrated planning algorithm, and the algorithms for ride-pooling, shared mobility, and private vehicles. Section \ref{Case Study} provides the experimental results of our model application obtained from applying our approach to a case study in the Rotterdam suburban area. Section \ref{Conclusions} concludes this study, discusses the advantages and limitations of our approach, and explores potential related avenues for future research.}

\section{Literature Review}\label{Literature Review}

{This section reviews various approaches and methodologies used in multi-modal transport planning.} {Past studies included in our review are selected based on a selection of recent papers that explore the integration of PT with SM. Keywords such as ``shared mobility", ``public transport integration", and ``preference-based optimization" guided our search in literature. Specifically, we prioritized studies that focused on optimization methods, user preferences, and system-level integration. Papers focusing solely on preferences or unimodal transport were filtered out, resulting in the identification of 12 highly relevant studies. These studies were classified into three main categories: Integrated Dial-A-Ride Problem (IDARP), Flexible Mobility on Demand (FMOD), and Choice-Driven Service Network Design (CSND). The IDARP category focuses on integrating ride-pooling services with PT, while FMOD explores flexible mobility services without fixed-line PT. In contrast, CSND considers only the integration of various PT modes. In this study, we bridge PT and SM, as well as supply and demand, to enable seamless multimodal trips.}

Some studies employ arc-based formulations in modeling integrated multi-modal transport planning. For example, \cite{hall2009integrated} introduce an arc-based formulation to solve the IDARP, which involves optimizing the scheduling of dial-a-ride requests alongside fixed-route PT services. It provides a solution for a small-scale instance of the problem, considering two vehicles, one PT line, and four ride requests, and simplifies the problem by overlooking waiting times during transfers and synchronization between transport modes. \cite{posada2017integrated} focus on arc-flow formulations that enable integration of demand-responsive vehicles and fixed-route transport. The study presents two distinct arc-flow models, differing in their approach to modeling transfer points, and explores ways to enhance these models through valid inequalities. \cite{varone2015multi} present a multi-modal routing problem, focusing on optimizing transport using a combination of ride-sharing and PT. \cite{huang2018multimodal} focus on optimizing multi-modal trips through the application of shortest path algorithms. They leverage two real-world datasets from the Swiss Federal Railways and a prominent European carpooling platform. They construct and query a network graph using Neo4j and use routing algorithms like the Dijkstra algorithm to find the shortest path.

{Given its inherent complexity, multi-modal planning usually involves several sequential steps. Nonetheless, it is important to highlight that the specific steps that are undertaken in multi-modal transport planning vary in the literature.} \cite{varone2015multi} propose a method that starts with a user's PT itinerary, selects potential drivers based on closeness estimation, computes driving paths, and determines the best solution for the earliest arrival time, with experiments conducted using real-world data from Swiss cities. \cite{stiglic2018enhancing} integrate ride-sharing and public transit through two phases: match identification and optimization. They separately identify matches for shared rides, transit, and park-and-ride transit. However, they do not consider matches with more than two riders.

{Recognizing the computational complexity involved in integrated planning of PT and ride-sharing, several studies propose (meta)heuristic algorithms as a means to expedite the computation process.} \cite{lee2019development} use a simulated annealing approach for bus routing between train stations, and also consider vehicle relocation based on demand. If optimal bus routes cannot be achieved with buses at specific stations, additional buses will be relocated to those stations. \cite{posada2020metaheuristic} use ALNS to optimize the routes of vehicles considering the timetables of PT. \cite{yu2021optimal} propose an angle-based clustering algorithm for multi-modal transport with PT and shuttles, considering the convenience of passengers in terms of driving time, travel time, and number of transfers. \cite{molenbruch2021analyzing} solve the static IDARP using the Large Neighborhood Search (LNS) and take into account the synchronization between shared vehicles and PT.

{A few studies delve into the real-time dimension of this domain.} For instance, \cite{ma2017demand} discuss a dynamic bi-/multi-modal vehicle dispatching and routing algorithm aimed at optimizing real-time ride-sharing (feeder) services alongside existing PT networks. The proposed algorithm focuses on shareability to group ride requests efficiently, allowing seamless bi-/multi-modal trips while considering factors like demand intensity, fleet size, waiting times, passengers' delays, and vehicle occupancy. {\cite{akyuz2023partial} study partial and complete re-planning of a multi-modal transport system under disruptions.}

{A related thread of studies concentrates on comparing the performance and features of PT and ride-sharing services rather than their integration.} \cite{fidanoglu2023integrating} solve a DARP and compare the solutions of shared autonomous vehicles and PT and find that service quality and energy efficiency can be improved using shared autonomous vehicles. However, they do not consider transfers between shared autonomous vehicles and PT. Considering PT and cars, \cite{esztergar2020activity} solve a traveling salesman problem for individual passengers, while they do not consider multiple passengers and passenger preferences.

{Passenger preferences are considered in the integrated planning of PT and SM in the study of \cite{azadeh2022choice}. {They investigate the integration of two types of bus services, i.e., Demand Responsive Transport (DRT) as a supplement to Fixed Line and Schedule (FLS), in less densely populated regions. DRT offers flexible, app-based transportation services with customizable routes and schedules, while FLS refers to public transport systems operating on predetermined routes and timetables, catering to regular travel patterns.} In their approach, passenger behavior is implicitly modeled through a discrete choice model, allowing for a nuanced understanding of user preferences and their impact on system design.} {However, their focus is on determining the locations of FLS and DRT stations, whereas our study concentrates on trip planning decisions.}

{\cite{atasoy2015concept} develop a model that offers passengers a menu of optimized travel options, allowing them to choose based on their preferences. The system's flexibility lies in allocating vehicles to different services throughout the day, enabling them to switch roles as needed. However, they do not consider transfers and multimodal transport.}

{Table \ref{table:literature} summarizes and compares the approaches reported in the literature as well as the approach undertaken in this study (last row). Most studies focus on IDARP \citep{hall2009integrated, varone2015multi,ma2017demand,posada2017integrated,huang2018multimodal,stiglic2018enhancing,posada2020metaheuristic,yu2021optimal,molenbruch2021analyzing,fidanoglu2023integrating}, {which consider ride-pooling services but do not consider micro-mobility services}, while this study addresses a Dynamic Multi-modal Trip Planning (DMTP) problem, integrating PT and various SM services. Only a few studies, such as \cite{varone2015multi} and \cite{ma2017demand}, have proposed dynamic planning approaches, and this study introduces a dynamic planning method based on the Rolling Horizon framework. Notably, prior studies often overlook passenger preferences in trip planning, whereas our approach accounts for the preferences of passengers. To the best of our knowledge, this study is the first to propose a dynamic multi-modal transport planning approach that explicitly incorporates passenger preferences.}

\begin{table}[h]

\begin{adjustbox}{max width=\textwidth}

 \begin{threeparttable}
\caption{Comparison between the proposed approach and existing approaches in the literature}
\label{table:literature}

\begin{tabular}{p{4cm} c c c c c c c}

\toprule

Article & Problem&Mode{$\dagger$}&Transfer$\star$&Dynamic&User preference&Methodology\\
\midrule

\cite{hall2009integrated}&IDARP&PT, RP&\checkmark&&&MIP\\
\cite{varone2015multi}&IDARP &{B, TR}, RP, W&\checkmark&RH&& SS\\
{\cite{atasoy2015concept}}&FMOD&T, RP, MB&&&\checkmark&MIP\\

\cite{ma2017demand}&IDARP&{TR}, RP&\checkmark&RH&&A*\\%
\cite{posada2017integrated}&IDARP&PT, RP, W&\checkmark&&&MIP\\
\cite{huang2018multimodal}&IDARP&{TR}, RP&\checkmark&&& Dijkstra and A*\\
\cite{stiglic2018enhancing}&IDARP&{TR}, RP, W&\checkmark&&&\\
\cite{lee2019development}&LMT&{TR}, RP&&&&RM\\%
\cite{posada2020metaheuristic}&IDARP&{B, TR}, RP&\checkmark&&&ALNS\\
\cite{yu2021optimal}&IDARP&{TR, S}, RP, W&\checkmark&&&AC\\%
\cite{molenbruch2021analyzing}&IDARP&PT, RP&\checkmark&&&LNS\\%
\cite{azadeh2022choice}&CSND&{B} (FLS, DRT)&\checkmark&&\checkmark&MIP and ALNS\\
\cite{fidanoglu2023integrating}&IDARP&{B}, RP, W&&&&NS\\%

\cmidrule(lr){1-7}
This article & DMTP  & {B, TR}, RP, SMM, W&\checkmark&RH &\checkmark&MIP and ALNS\\
\bottomrule
\end{tabular}

\begin{tablenotes}
      \footnotesize
      \item {$\dagger$: When the relevant study specifies the PT mode, it is listed explicitly. Otherwise, it is referred to as ``PT".}
      \item {$\star$: The ``Transfer" column indicates whether the relevant study considers the process of passengers transferring between different modes of transport. This includes factors such as waiting times during transfers and synchronization between transport modes.}
      \item IDARP: Integrated Dial-a-Ride Problem; FMOD: Flexible Mobility on Demand; LMT: Last Mile Transit; CSND: Choice-driven Service Network Design; DMTP: Dynamic Multimodal Trip Planning; {B: Bus; TR: Train;} PT: Public Transport; RP: Ride Pooling; {S: Subway;} SMM: Shared Micro-mobility; W: Walking; T: Taxi; MB: Mini Bus; FLS: Fixed Line and Schedule; DRT: Demand Responsive Transport; RH: Rolling Horizon.
      \item MIP: Mixed Integer Programming; SS: Sub-path
substitution; RM: Ride-matching algorithm; A*: A* algorithm; Dijkstra: Dijkstra algorithm; NS: Neighborhood Search; LNS: Large Neighborhood Search; ALNS: Adaptive LNS; AC: Angle-based Clustering algorithm.

    \end{tablenotes}
  \end{threeparttable}

  \end{adjustbox}
\end{table}

{The proposed framework integrates disaggregated demand parameters from an existing published survey \citep{montes2023shared} conducted in the study area to ensure compatibility. We model a network of PT and various SM services, enhancing the generalization capability of the framework. We introduce a choice-based optimization mixed-integer programming (MIP) model where demand assumptions remain independent of resource availability, and user preferences are not influenced by supply constraints. To address computational efficiency, we implement a Rolling Horizon approach for dynamic decision-making and employ a tailored heuristic algorithm to solve larger instances more efficiently. This uniquely designed framework distinguishes our approach from the current state of the art.}

\section{Problem description}
\label{Problem description}

In the following, we envisaged a preference-based multi-modal trip planner, which is designed to seamlessly provide PT and SM services (e.g. shared (e-)bikes, or scooters) to passengers. {Walking to/from PT stations and SM vehicles is also considered.} Table \ref{Notation} lists the notations used in this paper.

{In this study, we define SM as mobility solutions in which users share a fleet of vehicles, either concurrently (e.g., ride-pooling services where multiple passengers share a vehicle simultaneously) or sequentially (e.g., micro-mobility services such as shared bikes and e-scooters that are used by different individuals at different times). We treat micro-mobility—a subset of shared mobility—as a distinct category, encompassing lightweight, small-footprint vehicles. In our study, micro-mobility focuses on shared bikes and scooters as they are commonly used in urban environments for either first- and last-mile connectivity or door-to-door services. The term scooter is broader than moped and includes mopeds as a subset. Here, we use the term scooter because we are not referring to a specific type but rather emphasizing its key characteristic of being typically free-floating. Our proposed approach accommodates both dock-based bikes and free-floating shared scooters, making it adaptable to both dock-based and free-floating micro-mobility systems.}

Each passenger initiates a request ($r$) through the platform, with request details encompassing origin ($p_r$), destination ($d_r$), time window ($[a_r, b_r]$), and the number of passengers ($q_r$) to be transported. {Passengers may voluntarily register their information, such as gender, age, familiarity with shared mobility, and education level. If a passenger belongs to a specific segment, the corresponding preference data will be applied. Preference data can be obtained through a combination of data collection techniques, including surveys and questionnaires to capture user preferences, data from mobility apps and public transport platforms tracking real-time behavior, and open government datasets on transport usage. However, determining preference parameters falls outside the scope of this study.} {The utility parameters are entirely exogenous to the decision-making process of a trip planner from the MaaS platform perspective. The trip planner considers inventory and fleet availability, using predefined user preferences to design multimodal trips for each request. We do not revisit or adjust individual utility parameters after the trip assignment. Instead, the system optimizes mode assignment based on availability while adhering to the initially defined preferences.}

After a passenger initiates a request, the platform recommends an alternative from the available options. {We consider both PT and SM services to offer door-to-door services.} The platform is able to provide several types of services, including:
 \begin{enumerate}
   \item PT and sharing micro-mobility,
\item PT and ride-pooling vehicles,
\item PT and walking,
\item sharing micro-mobility,
\item and private vehicles/walking as opt-out options.
 \end{enumerate}

{The category ``PT and sharing micro-mobility" encompasses both PT with shared bikes and PT with shared scooters. When both PT and SM are used, SM is utilized for the first-mile and last-mile connection, and public transport serves as the primary leg of the journey. If PT is not used, a single SM mode can provide a complete door-to-door connection. While combinations like shared scooters and shared bikes are not included, such combinations could be explored in future studies to enhance flexibility.} {In this study, PT refers to fixed-line public transport, which represents the most common form of public transit. The selection of SM solutions is motivated by the need for a generalizable framework. SM can be interpreted either as a fleet assigned per passenger (e.g., shared bikes and scooters) or as vehicles shared among multiple passengers (e.g., ride-pooling vehicles), with the distinction lying in the definition of capacity. This approach ensures that our model remains adaptable to various shared mobility configurations.} These diverse service options cater for a diverse pool of passengers with various preferences and needs.

{When a passenger initiates a request, the backend generates alternatives that not only respect operational constraints, such as micromobility vehicle occupancy, but also aim to maximize overall social welfare. Simply considering constraints is insufficient, as failing to account for social welfare could negatively impact the utility of other passengers, especially in ride-pooling scenarios where their routes might change. {By considering both constraints and social welfare, the backend generates multiple alternatives and selects the optimal option for the passenger considering their preferences.}}

The platform operates dynamically. It receives passenger requests, generates alternative solutions using a range of transport modes, and then either assigns a ride-pooling vehicle or directs the passenger to a micro-mobility or public transport service location considering passengers' preferences and profiles. In the process of integrating PT with ride pooling or micro-mobility services, the PT services are initially matched with the request. {Subsequently, based on the availability of PT services, choices are made between shared mobility options such as ride-pooling and micro-mobility services.}

\begin{xltabular}{\textwidth}{l X}
\caption{Notation.}\label{Notation}\\

\hline \multicolumn{2}{l}{\textbf{Sets:}}\\
\endfirsthead

  $W$ & Set of transport modes, $w \in W$.\\
  $R$ & Set of requests, $r \in R$. $R^{t} \subseteq R$ denotes the set of active requests at time $t$, which includes those that have not yet been scheduled or have not yet reached their intended destinations. $R_{\textrm{finish}} \subseteq R$ represents requests that already reached their destinations. \\
  $N$ & Set of locations, $i, j \in N$. $O/P/D/T \subseteq N$, set of depots/origins/destinations/transfer locations. $C^t \subseteq N$, set of current locations of vehicles at time step $t$. $N_{w} \subseteq N$ represents the set of locations for transport mode $w$.\\

  $K$ & Set of vehicles, $k \in K$. $K_{\textrm{\textrm{PT}}}\subseteq K$, $K_{\textrm{\textrm{ride}}}\subseteq K$, and $K_{\textrm{micro}}\subseteq K$ represent sets of PT services, ride-sharing vehicles, and micromobility vehicles, respectively. $K_{\textrm{private}}$ and $K_{\textrm{private}}^r$ are sets of private vehicles and private vehicles owned by {passengers associated with request $r$}, respectively. \\

  $E$ & Set of routes, $(i,j) \in E$.\\

  $A$ & {Set of alternatives. $A_r$ represents the set of alternatives for request $r$.}\\
  $H$ & Set of schedules of finished requests $R_{\textrm{finish}}$, $(i,j,k,r) \in H$.\\

 \hline
\multicolumn {2}{l}{\textbf{Parameters:}}\\
 $e_k$ & Capacity of vehicle $k$.\\

 $q_r$ & Number of passengers associated with request $r$. \\
 $\tau_{ij}^k$ & Travel time [minute] on the shortest path between locations $i$ and $j$ for vehicle $k$.\\

 $[a_r, b_r]$ & Earliest starting time and latest ending time for request $r$.\\

 $v_k$ & Vehicle $k$'s speed [km/h]\\
 $v_r$ & Walking speed [km/h] {shared by all passengers associated with} request $r$.\\
 $d_{ij}^k$ & The distance [km] from location $i$ to location $j$ for vehicle $k$.\\

 $c_k^{n}$ & The transport cost [euro] per minute or kilometer for vehicle $k \in K$ is denoted as $c_k^1$/$c_k^{1’}$. The cost per hour of waiting time is denoted as $c_k^2$. \\
$M$ & A large enough positive number.\\
$ASC_{\text{w}}$ &Alternative specific constant of using transport mode(s) $m$.\\
$\beta$ &Parameters in travellers' utility functions. \\
 $t_{\textrm{search}}$ & Time for searching parking space.\\
 $c_{\textrm{park}}$ & Parking cost.\\

 \hline
  \multicolumn {2}{l}{\textbf{Variables:}}\\
 $x_{ij}^{k}$ &Binary variable indicating whether vehicle $k$ uses the edge $(i,j)$; 1 if it does, 0 otherwise.\\
 $y_{ij}^{kr}$ & Binary variable indicating whether request $r$ transported by vehicle $k$ between locations $(i,j)$; 1 if it does, 0 otherwise. \\

 $s_{ir}^{kl}$ & Binary variable indicating whether request $r$ is transferred from vehicle $k$ to vehicle $l$ at location $i$; 1 if it is, 0 otherwise. \\
 $g_{r}$ & Binary variable indicating whether request $r$ is transferred during transport; 1 if it is, 0 otherwise.\\
 ${x'}_r^k$ & Binary variable indicating whether request $r$ is transported by vehicle $k$; 1 if it is, 0 otherwise.\\
 $t_i^{kr}/{t'}_i^{kr} /\overline{t}_i^{kr}$ & The arrival time, beginning time, and departure time of request $r$ served by vehicle $k$ at location $i$.\\
 $t_i^k/{t'}_i^k/\overline{t}_i^k$ & The arrival time, beginning time, and departure time of vehicle $k$ at location $i$. \\
 $t_{ki}^{\textrm{wait}}$ & The duration vehicle $k$ waits at location $i$. \\

 $t^k$ & Micro-mobility vehicle $k$'s last operation time.\\
 $l^k$ & Micro-mobility vehicle $k$'s last location.\\
 $V_r^a$ & Deterministic utility of request $r$ and alternative $a$.\\

 $vt_w$ & In-vehicle time when using vehicles of transport mode $w$ or walking time when only walking.\\
 $cost_w$ & Cost of using transport mode(s) $w$.\\
$walk_{w}$ & Walking time using transport mode $w$.\\
$wt_{w}$ & Waiting time using transport mode $w$.\\

 \hline

\end{xltabular}

{PT adheres to real-world schedules that are predefined, meaning that we do not account for unexpected events and the re-planning of PT schedules. PT vehicles stop at each station for a certain amount of time to allow passengers to board and alight. Passengers have the option to transfer from SM to PT or from PT to SM. If an SM vehicle cannot reach the PT station, walking between the location of an SM vehicle and a PT station, or vice versa, is also allowed. Additionally, passengers can choose to walk from their point of origin to the PT station or from the PT station to their destination, i.e., walking may be used as an access and/or egress mode.}

In the case of ride-pooling, various types of vehicles, such as vans, taxis, or mini-buses, can be used with capacity $e_k$ and speed $v_k$.
{In our problem formulation, ride-pooling vehicles have the flexibility of either adhering to predefined routes or being flexible according to the demand. In the case of predefined routes, they adhere to schedules with time windows at each station. However, in the flexible mode, these vehicles have the freedom to traverse any route necessary to pick up and deliver passengers.} Ride-pooling services involve a fleet of vehicles initially stationed at depots (possibly corresponding to drivers' home locations in the case of ride-sourcing platforms). Each vehicle can cater to multiple requests simultaneously, provided its capacity is not exceeded. For example, a vehicle may pick up a passenger at location A, proceed to location B for another pickup, and eventually drop these passengers at locations C and D, respectively. Requests that need routing at time $t$ result in updated vehicle routes, with new requests inserted en-route into the existing schedules. In cases of insufficient capacity, some requests may remain unserved. When inserting new requests, the schedules of planned requests can be modified as long as the planned requests are not delayed due to these changes, and the order of planned requests remains unchanged.

Our model assumes no competition among service providers, and we operate under the assumption that drivers of ride-pooling services will accept tasks assigned by the platform. Additionally, the battery status of micro-mobility vehicles is not considered. We assume that transfers are limited to interchanges between SM and PT and within PT, i.e., no transfers between successive SM legs.

\section{{Rolling horizon approach for dynamic requests}}
\label{Rolling horizon approach}

\begin{figure}[H]
    \centering
    \includegraphics[scale=0.5]{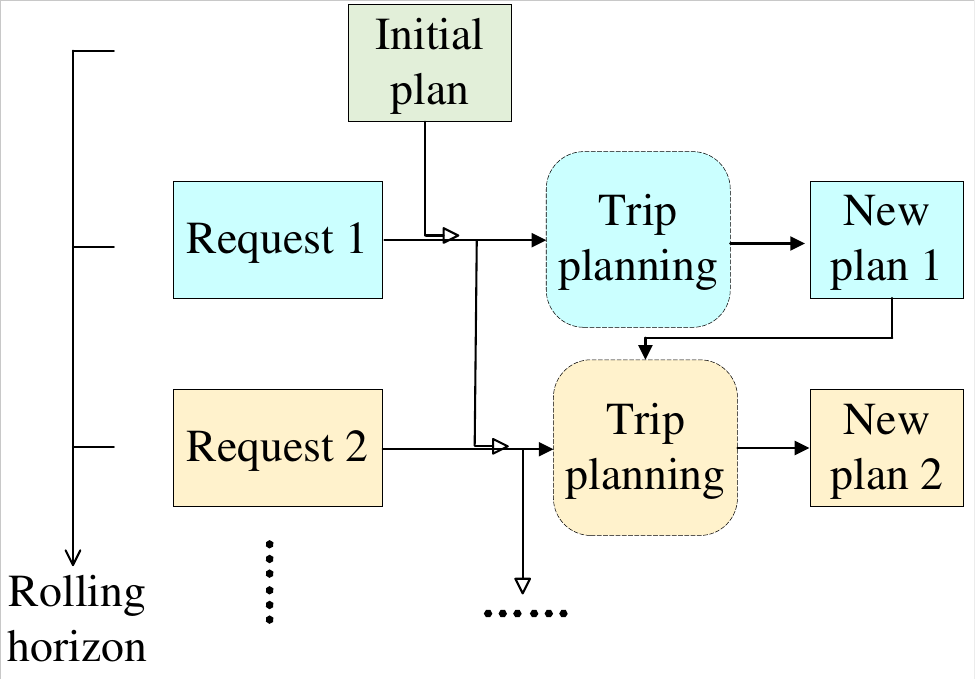}
    \caption{Dynamic multimodal trip planning based on rolling horizon.}
    \label{fig:Dynamic multimodal trip planning based on rolling horizon}
\end{figure}

In the real world, passengers request services at various points in time, ride-sharing drivers continuously announce their trips, and micro-mobility service providers continuously update the real-time positions of micro-vehicles \citep{stiglic2018enhancing}. {To effectively manage these dynamics, we implement a rolling horizon mechanism, which means that vehicles are planned dynamically when a new request is received. The proposed multimodal trip planner provides the best trip option to a single request at each time step, and we do not process batch requests within a single time step.} The best solution will be considered as the initial {plan} in the next time step, as shown in Figure \ref{fig:Dynamic multimodal trip planning based on rolling horizon}.

{In our rolling horizon approach, advanced reservations are applied to SM services, while PT is assumed to have unlimited capacity, negating the need for reservations. For SM services, the system assigns specific time slots for vehicles based on passenger requests. For micro-mobility, once a time slot is assigned to a passenger, it is exclusively reserved and cannot be used by others. In contrast, ride-pooling vehicles can accommodate additional passengers as long as there is remaining capacity. The vehicle remains available for scheduling outside the reserved time slot. For instance, if a shared bike is recommended to and taken by passenger A, it remains unavailable for passenger B while passenger A is using it. This design ensures that, at every time step, the system identifies the optimal service assignment for the incoming request. By structuring the system in this way, we guarantee that there is no inherent conflict between individual utility and system optimum since the optimization occurs at the level of each individual request without revisiting past decisions.}

\section{{Mathematical model for preference-based multi-modal trip planning}}
\label{Mathematical model for preference-based multimodal trip planning}
{The platform acts as a responsive mediator between passengers and available services. Initially, all passengers using this platform are asked to provide their profile information. However, passengers may choose not to provide this information. In such cases, the platform will consider the passenger to not have specific preferences which are known a-priori. Subsequently, when the passenger requires services, his request encompasses essential details, including origin, destination, desired departure time, and the latest acceptable arrival time. {This approach is detailed in Figure \ref{fig:multimodal trip planning model}, which illustrates a multi-modal trip planning framework that incorporates heterogeneous passenger preferences into the planning process. It begins with journey planning, where diverse passenger preferences, influenced by socio-demographic characteristics and travel objectives, are modeled. Requests from passengers are processed through a multimodal transport platform, which matches them with available services. The preference-based trip planning integrates passenger utility functions into an optimization model to identify the most suitable transport alternatives, balancing factors such as cost, time, and satisfaction. The output is a transport plan that considers preferences across multiple passengers, generating optimized routes and connections tailored to their individual needs.}

\begin{figure}[h]
    \centering
    \includegraphics[scale=0.55]{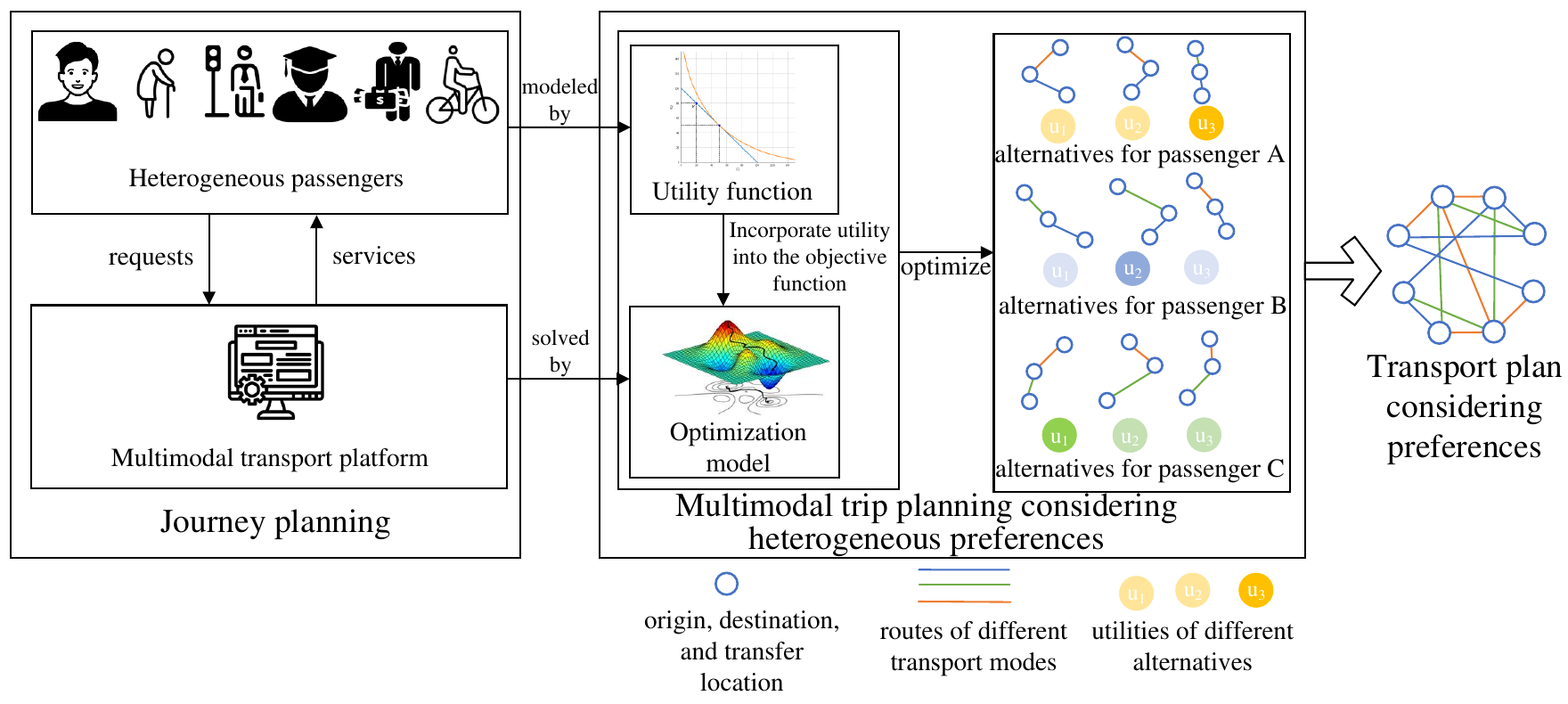}
    \caption{Multi-modal trip planning model.}
    \label{fig:multimodal trip planning model}
\end{figure}

This section first introduces how the utility is calculated in Section \ref{Utility calculation}, and then proposes a mathematical model for multi-modal trip planning in Section \ref{Mathematical model for multimodal trip planning}.}

\subsection{{Heterogeneous passenger mobility preferences}}
\label{Utility calculation}
{{Preference-based optimization can be categorized in three ways: (a) assuming all passengers have homogeneous preferences, (b) recognizing that different passenger segments have heterogeneous preferences, or (c) considering individual preferences in detail. In practice, it is difficult to gather precise individual preferences, while option (a) is overly simplistic and fails to reflect real-world conditions. Therefore, this study adopts option (b) as a balanced approach. Conducting surveys to gather preferences is a common method in the literature on preference analysis; however, this is beyond the scope of the current study.} To incorporate passengers' preferences, we {adopt} choice model estimations based on a stated preference survey reported in \cite{montes2023shared}}. \cite{montes2023shared} study the multi-modal transport in Rotterdam, focusing on the preferences and choices associated with both PT usage and the subsequent egress modes, which encompass shared bikes, shared scooters, and walking. {Our assumptions regarding preferences are as follows:}

\begin{enumerate}
    \item 	{Homogeneity within each group:
We assume that individuals exhibit similar preferences within each user group. Socio-demographic characteristics such as age, frequency of public transport use, and previous experience with shared micromobility influence mode choice.}
\item {Attributes: Attributes considered include travel cost, in-vehicle time, waiting time, and walking time. The study differentiates between PT attributes and SM attributes, with separate cost sensitivity for each. Unobserved attributes such as comfort and convenience are captured through ASC.}
\item {Alternatives: We evaluate multiple transport alternatives, as detailed in Section 3, with choice sets tailored to passengers based on service availability.}
\item {Choice mechanism: The alternative with the highest utility will be chosen.}
\end{enumerate}

{These assumptions are motivated by the need to accurately capture real-world decision-making in multimodal transport. Preferences are modeled based on fundamental trade-offs between travel time and cost, as these are widely recognized as primary determinants of mode choice. The inclusion of different passenger groups is necessary to account for variations in behavior, particularly differences in sensitivity to cost, time, and familiarity with shared mobility. In addition, we assume static user preferences, meaning that individual choices are based on pre-determined parameters without real-time adjustments. However, the dynamic aspect is incorporated on the supply side, reflecting variations in resource availability, such as the presence of SM fleets at PT stations. This structure allows for a generalizable model while maintaining computational tractability.}

{\cite{montes2023shared} estimate Multinomial Logit (MNL) models with interaction effects of passengers' profile information, which includes socio-demographic characteristics and transport-related information. The considered socio-demographic characteristics include gender, age, education, and income. The transport-related information includes frequency of use of PT, familiarity with SM, and previous use of SM.} {To obtain information on passengers' sociodemographic data, the platform can request such details from users upon engagement. In some countries, such as the Netherlands, sociodemographic information is also publicly accessible \citep{NL_socio2,NL_socio}. If specific data is unavailable, the platform can apply general preference parameters.} {Preferences are defined based on the market share of different alternatives for each user group, reflecting their individual needs and priorities regarding different attributes of the transportation options, such as cost and travel time. The survey conducted by \cite{montes2023shared} is in the same area as this study, and we utilize the results—specifically, the estimated preference parameters for different passenger profiles—as prior information. The profiles are illustrated in detail in Section \ref{Numerical experiments and managerial insights}.}

{We consider passengers’ preferences regarding the cost and travel time attributes of trip alternatives, along with intrinsic preferences for factors such as comfort and convenience. The deterministic utility of each request $r$ and alternative $a$ is calculated based on Equation \ref{cost1}:}

\begin{dmath}
V_r^a = ASC + \beta^{\textrm{cost}}cost + \beta^{\textrm{time}}time \label{cost1}
\end{dmath}
{where ASC denotes the Alternative Specific Constant, capturing intrinsic preferences for specific modes.}

For each request, the platform identifies the optimal alternative with the highest utility, taking into account the overall social welfare of both current and previous passengers. The utility of the alternative proposed to the passenger is calculated using Equation \ref{utility_from_system} \citep{montes2021studying, montes2023shared}:

\begin{dmath}
    {V_r =  \sum_{k \in K} ASC_r^{k} {x'}_r^k + g_{r} \sum_{(i,j) \in E} ( \sum_{k \in K_{\textrm{PT}}}  y_{ij}^{kr} (\beta_r^{PTwait} t_{ki}^{\textrm{wait}} + \beta_r^{\textrm{mainTime}} \tau_{ij}^k + \beta_r^{\textrm{mainCost}} cost_{ij}^k)} \\{+ g_{r} \sum_{k \in K \setminus K_{\textrm{PT}}} y_{ij}^{kr} (\beta_r^{\textrm{subTime}} (t_{ki}^{\textrm{wait}} + \tau_{ij}^k) + \beta_r^{\textrm{subCost}} cost_{ij}^k))}\\ {+ (1-g_{r}) \sum_{(i,j) \in E} (\sum_{k \in K_{\textrm{walk}}}   y_{ij}^{kr} \beta_r^{\textrm{walk}} \tau_{ij}^k  + \sum_{k \in K \setminus
K_{\textrm{walk}}}
 y_{ij}^{kr} (\beta_r^{\textrm{mainTime}} \tau_{ij}^k +
 \beta_r^{\textrm{mainCost}} cost_{ij}^k) )}
 \label{utility_from_system}
\end{dmath}

{where $\beta_r^{PTwait}$, $\beta_r^{\textrm{mainTime}}$, $\beta_r^{\textrm{mainCost}}$, $\beta_r^{\textrm{subTime}}$, and $\beta_r^{\textrm{subCost}}$ represent $\beta$ parameters for waiting time at PT stations, travel time using main transport mode, travel cost using main transport mode, travel time for access and egress legs, and travel cost for access and egress legs, respectively. The main transport mode means the primary mode of transport, whether it is PT in a multimodal trip or shared micro-mobility in an unimodal trip. The decision variables in Equation \ref{utility_from_system}, including ${x’}_r^k$ and $y_{ij}^{kr}$, are determined by the optimization model in Section \ref{Mathematical model for multimodal trip planning} to identify the alternative with the highest utility.}

The cost is a function of time and distance, as calculated by Equation \ref{cost}. The considered cost includes the fixed cost ($cost_{fixed}^k$), cost for travel time ($\tau_{ij}^k * cost_{time}^k$), and cost for travel distance ($d_{ij}^k * cost_{distance}^k$):
\begin{dmath}
cost_{ij}^k = cost_{fixed}^k + \tau_{ij}^k * cost_{time}^k + d_{ij}^k * cost_{distance}^k \label{cost}
\end{dmath}
{where $\tau_{ij}^k$ and $d_{ij}^k$ represent the actual travel time and distance using vehicle $k$, respectively.}

\subsection{Mathematical model for multi-modal trip planning}
\label{Mathematical model for multimodal trip planning}

{The platform provides a trip recommendation as soon as it receives a request $r_t$ at time $t$. For requests ${R^t}\setminus r_t$ that have already been assigned services in previous time steps, their designated services remain unchanged, meaning only the current request $r_t$ needs to be planned. The objective of the platform is to maximize the social welfare of all passengers $R^t$ who have not arrived at their destinations:}

\begin{dmath}
\max {F(t) = \sum_{r \in {R^t}} V_r = \sum_{r \in {R^t}\setminus r_t} V_r + \sum_{a \in A_{r_t}} {x'}_{r_t}^a {V_{r_t}^a} } \label{obj2}
\end{dmath}

{where ${x'}_{r_t}^a = 1$ indicates that the alternative $a$ is selected for request $r$, {$V_r$ represents the utility of request $r$, and $V_{r_t}^a$ represents the utility of alternative $a$ for request $r_t$.}} {The system serves passengers sequentially in a dynamic setting and provides an optimal solution for the individual passenger while accounting for resource constraints imposed by services previously assigned to other passengers.}

Constraints \eqref{served_requests} to \eqref{wait time} are constraints for all transport modes. Constraints \eqref{served_requests} make sure that the number of served requests reaches a predefined level, which can be adjusted by parameter $\varepsilon$. The passenger $r$ is considered served when departing from the origin $p_r$ with a vehicle $k$ ($y_{p_rj}^{kr}$ = 1). Constraints \eqref{begin request} ensure that if a passenger starts their trip, they must also complete it by reaching their destination.
  \begin{align}
  &\sum_{r \in {R^t}}\sum_{k \in K}\sum_{j \in N} y_{p_rj}^{kr} \geqslant \varepsilon |{R^t}| & \label{served_requests}\\
&\sum_{k \in K}\sum_{i \in N} y_{id_r}^{kr} \geqslant \sum_{k \in K}\sum_{j \in N} y_{p_rj}^{kr} \quad \forall r \in {R^t} & \label{begin request}
\end{align}

Constraints \eqref{begin} enforce that each vehicle may initiate at most one route from its initial location; Constraints \eqref{end} ensure that if a vehicle is used, it ends the route at its designated location, except for micro-mobility and walking. Constraints \eqref{flow conservation of vehicle} represent vehicle flow conservation. Constraints \eqref{vehicle flow cover request flow} link $y_{ij}^{kr}$ and $x_{ij}^{k}$ variables in order to guarantee that for a request to be transported by a vehicle, that vehicle needs to traverse the associated route.
\begin{align}
&\sum_{j \in N} x_{o_kj}^{k}  \leqslant 1 \quad  \forall k \in K\setminus K_{\textrm{walk}}& \label{begin}\\
&\sum_{j \in N} x_{o_kj}^{k} = \sum_{j \in N} x_{jo'_k}^{k} \quad \forall k \in K \setminus K_{\textrm{walk\&micro}} & \label{end}\\
&\sum_{j \in N} x_{ij}^{k} - \sum_{j \in N} x_{ji}^{k}=0 \quad \forall k \in K\setminus K_{\textrm{walk}}, \text{ } \forall i \in N\setminus{o_k,o'_k}& \label{flow conservation of vehicle}\\
&y_{ij}^{kr} \leqslant x_{ij}^{k} \quad \forall (i,j) \in E,\text{ } \forall k \in K,\text{ } \forall r \in {R^t}& \label{vehicle flow cover request flow}
\end{align}

Constraints \eqref{transshipment1} ensure that there is only one transfer for one request at a given transfer location. Constraints \eqref{no_transshipment_when_no_PT}, \eqref{transshipment3}, and \eqref{transshipment4} forbid transfers without PT services and walking, transfers between micro-mobility vehicles, and transfers between the same vehicle $k$, respectively.
\begin{align}
&\sum_{j \in N} y_{ji}^{kr} + \sum_{j \in N} y_{ij}^{lr} \leqslant s_{ir}^{kl} +1 \quad \forall r \in {R^t}, \text{ } \forall i \in T,  \text{ } \forall k,l \in K& \label{transshipment1}\\
&s_{ir}^{kl} = 0 \quad \forall r \in {R^t}, \text{ } \forall i \in T, \text{ } \forall k,l \in K\setminus K_{\textrm{PT\&walk}}& \label{no_transshipment_when_no_PT}\\
&s_{ir}^{kl} = 0 \quad \forall r \in {R^t}, \text{ } \forall i \in T, \text{ } \forall k \in K_{\textrm{micro}}, \text{ } \forall l \in K_{\textrm{micro}} & \label{transshipment3}\\
&s_{ir}^{kk} = 0 \quad \forall r \in {R^t}, \text{ } \forall i \in T, \text{ } \forall k \in K & \label{transshipment4}
\end{align}

Constraints \eqref{capacity} are the capacity constraints.
\begin{align}
&\sum_{r \in {R^t}}q_r y_{ij}^{kr} \leqslant e_k x_{ij}^{k} \quad \forall (i,j) \in E,\text{ } \forall k \in K& \label{capacity}
\end{align}

Constraints \eqref{right_T} ensure the transfer occurs in the PT stations that have the available facilities for transfers, such as SM hubs.
\begin{align}
&s_{ir}^{kl} = 0 \quad \forall k \in K_{w_1}, \text{ } \forall l \in  K_{w_2},\text{ } \forall i \in T\setminus T_{w_1}^{w_2}, \text{ } \forall r \in {R^t}, \text{ } \forall w_1, w_2 \in W & \label{right_T}
\end{align}

Constraints \eqref{beginning time} ensure that the beginning time occurs after the arrival time of passengers. Constraints \eqref{departure time} ensure that the departure time equals beginning time plus duration ${t''}_i^{kr}$. The time ${t''}_i^{kr}$ varies according to the specific transport mode. In the context of public transport, it means dwelling time, whereas in the case of micro-mobility, it corresponds to the time required for unlocking and locking the vehicle. Constraints \eqref{vehicle departure time} ensure that departures occur only after all passengers have boarded or disembarked from the vehicle. Constraints \eqref{vehicle arrives earlier than request arrival time} ensure that the request and vehicle's arrival time are the same. Constraints \eqref{vehicle last beginning time} define the beginning time.

\begin{align}
& t_i^{kr} \leqslant  {t'}_i^{kr} \quad \forall i \in N, \text{ } \forall k \in K, \text{ } \forall r \in {R^t}& \label{beginning time}\\
& {t'}_i^{kr} + {t''}_i^{kr} \sum_{j \in N} y_{ij}^{kr} \leqslant  \overline{t}_i^{kr} \quad \forall i \in N, \text{ } \forall k \in K, \text{ } \forall r \in {R^t}& \label{departure time}\\
&\overline{t}_i^k \geqslant \overline{t}_i^{kr}  \quad \forall i \in N, \text{ } \forall k \in K, \text{ } \forall r \in {R^t}& \label{vehicle departure time}\\
&t_i^k \leqslant t_i^{kr} \quad \forall i \in N,\text{ } \forall k \in K,\text{ } \forall r \in {R^t}& \label{vehicle arrives earlier than request arrival time}\\
&{t'}_i^k \geqslant {t'}_i^{kr} \quad \forall i \in N,\text{ } \forall k \in K,\text{ } \forall r \in {R^t}& \label{vehicle last beginning time}
\end{align}

Constraints \eqref{time on arc 1} and \eqref{time on arc 2} ensure that the travel time is consistent with the distance traveled and speed. Constraints \eqref{time window for pick up} take care of the time window for the origin and destination.

\begin{align}
&\overline{t}_i^k+\tau_{ij}^k-t_j^k \leqslant M(1-x_{ij}^k) \quad \forall (i,j) \in E, \text{ } \forall k \in K\setminus K_{\textrm{fix}} & \label{time on arc 1}\\
&\overline{t}_i^k+\tau_{ij}^k-t_j^k \geqslant -M(1-x_{ij}^k) \quad \forall (i,j) \in E, \text{ } \forall k \in K\setminus K_{\textrm{fix}} & \label{time on arc 2}\\
& {t'}_{p_r}^{kr} \geqslant a_r y_{ij}^{kr}, \text{ } \overline{t}_{d_r}^{kr} \leqslant b_r \quad \forall (i,j) \in E, \forall r \in {R^t}, \text{ } \forall k \in K & \label{time window for pick up}
\end{align}

Constraints \eqref{two vehicle's time on transshipment terminal} are time constraints for transfers. If a transfer occurs, the boarding time ${t'}_i^{lr}$ for vehicle $l$ must be after the alighting time $\overline{t}_i^{kr}$ from vehicle $k$. Constraints \eqref{wait time} calculate passenger waiting time.

\begin{align}
&\overline{t}_i^{kr}-{t'}_i^{lr} \leqslant  M(1-s_{ir}^{kl}) \quad \forall r \in {R^t}, \text{ } \forall i \in T,  \text{ } \forall k,l \in K,  \text{ } k \neq l & \label{two vehicle's time on transshipment terminal}\\
&t_{ki}^{\textrm{wait}} \geqslant {t'}_i^k - t_i^k \quad \forall i \in N,\text{ } \forall k \in K& \label{wait time}
 \end{align}

The constraints represented by \eqref{PT_must_run} guarantee the continuous operation of public transport vehicles, even in the absence of passengers.
\begin{align}
& x_{ij}^{k} = 1 \quad \forall (i,j,k) \in PT& \label{PT_must_run}
\end{align}

Constraints \eqref{flow conservation of request at regular terminal}-\eqref{flow conservation of request at T when not transshipment2} are dedicated to PT and ride-pooling services. Constraints \eqref{flow conservation of request at regular terminal} to \eqref{flow conservation of request at T when not transshipment2} represent request flow conservation. Specifically, Constraints \eqref{flow conservation of request at regular terminal} and \eqref{flow conservation of request at transshipment terminal} apply to regular and transfer locations, respectively. If request $r$ is not transferred at location $i \in T$ but vehicle $k$ passes location $i$ for other requests, Constraints \eqref{flow conservation of request at regular terminal} do not apply for request $r$. Therefore, additional flow conservation constraints (Constraints \eqref{flow conservation of request at T when not transshipment1} and \eqref{flow conservation of request at T when not transshipment2}) are included.

\begin{align}
\begin{split}
    & \sum_{j \in N_{\textrm{ride\&public}}} y_{ij}^{kr} - \sum_{j \in N_{\textrm{ride\&public}}} y_{ji}^{kr} = 0 \\
    &\quad \hspace{6cm} \forall k \in K_{\textrm{ride\&public}}, \text{ } \forall r \in {R^t}, \text{ } \forall i \in N_{\textrm{ride\&public}}\setminus T, p_r, d_r
    \end{split}\label{flow conservation of request at regular terminal}\\
\begin{split}
&\sum_{k \in K_{\textrm{ride\&public}}}\sum_{j  \in N_{\textrm{ride\&public}}} y_{ij}^{kr} - \sum_{k \in K_{\textrm{ride\&public}}}\sum_{j \in N_{\textrm{ride\&public}}} y_{ji}^{kr} = 0 \\
&\quad  \hspace{6cm} \forall k \in K_{\textrm{ride\&public}}, \forall r \in {R^t}, \forall i \in T \setminus p_r, d_r
\end{split}
\label{flow conservation of request at transshipment terminal}\\
\begin{split}
&\sum_{j \in N_{\textrm{ride\&public}}} y_{ij}^{kr} - \sum_{j \in N_{\textrm{ride\&public}}} y_{ji}^{kr} \leqslant \sum_{l \in K_{\textrm{ride\&public}}} s_{ir}^{lk}\\
&\quad \hspace{6cm} \forall k, l \in K_{\textrm{ride\&public}},\text{ } \forall r \in {R^t},\text{ } \forall i \in T \setminus p_r, d_r \end{split}\label{flow conservation of request at T when not transshipment1}\\
\begin{split}
&\sum_{j \in N_{\textrm{ride\&public}}} y_{ji}^{kr} - \sum_{j \in N_{\textrm{ride\&public}}} y_{ij}^{kr} \leqslant \sum_{l \in K_{\textrm{ride\&public}}} s_{ir}^{kl}\\
&\quad \hspace{6cm} \forall k, l \in K_{\textrm{ride\&public}},\text{ } \forall r \in {R^t},\text{ } \forall i \in T \setminus p_r, d_r \end{split}\label{flow conservation of request at T when not transshipment2}
\end{align}

Constraints \eqref{fixed route} ensure that vehicles operate along predefined routes or park at designated locations. Routes of PT are predefined, and ride-pooling vehicles may also operate along predetermined routes. Dock-based micro-mobility vehicles can only be parked in predefined locations, therefore, Constraints \eqref{fixed route} also apply to these vehicles. {This constraint is critical for maintaining operational consistency and ensuring that these services are available as planned for passengers. Without this constraint, the model could mistakenly allow vehicles (e.g. buses and trains) to deviate from their predefined routes.} Constraints \eqref{time window for fixed terminal} ensure that vehicles follow time windows for predefined stations, and they are not applied to micro-mobility vehicles.
\begin{align}
    &x_{ij}^{k}=0 \quad \forall k \in K_{\textrm{fix}},  \text{ } \forall (i,j) \in E\setminus {E}_{\textrm{fix}}^k& \label{fixed route}\\
    & t_i^{kr} \geqslant a_i^k y_{ij}^{kr},\text{ } \overline{t}_i^{kr} \leqslant b_i^k + M(1-y_{ij}^{kr})\quad \forall (i,j) \in E,\text{ } \forall r \in {R^t}, \text{ } \forall k \in K_{\textrm{fix}}\setminus  K_{\textrm{micro}} & \label{time window for fixed terminal}
\end{align}

Constraints \eqref{request flow cover vehicle flow} ensure that a micro-mobility is not able to move independently, as it is only allowed to be relocated by a passenger. {This reflects the real-world operational limitation where micro-mobility lacks autonomous movement capabilities and relies on users for repositioning. By enforcing this rule, the model accurately represents the practical constraints of micro-mobility systems and avoids unrealistic assumptions, such as vehicles relocating themselves.}
\begin{align}
&x_{ij}^{k} \leqslant y_{ij}^{kr} \quad \forall (i,j) \in {E}_{\textrm{micro}},\text{ } \forall k \in K_{\textrm{micro}},\text{ } \forall r \in {R^t}& \label{request flow cover vehicle flow}
\end{align}

{Constraints \eqref{finished_requests} are designed to ensure that previously assigned services for requests remain unchanged}.
\begin{align}
& y_{ij}^{kr} = 1 \quad \forall (i,j,k,r) \in H& \label{finished_requests}
\end{align}

Since a passenger can be transferred more than once, Constraints \eqref{transfer_or_not1} and \eqref{transfer_or_not2} set the binary variable $g_{r}$, which determines whether request $r$ is transferred or not. Constraints \eqref{k_transport_r} define the binary variable ${x'}_r^k$, which equals to 1 if a vehicle $k$ transport request $r$.

\begin{align}
& g_{r} \leqslant s_{ir}^{kl} \quad \forall k, l \in K \setminus K_{\textrm{walk}},\text{ } \forall r \in {R^t},\text{ } \forall i \in T& \label{transfer_or_not1}\\
& g_{r} M \geqslant s_{ir}^{kl} \quad \forall k, l \in K \setminus K_{\textrm{walk}},\text{ } \forall r \in {R^t},\text{ } \forall i \in T& \label{transfer_or_not2}\\
& {x'}_r^k \geqslant x_{ij}^k \quad \forall k \in K,\text{ } \forall r \in {R^t},\text{ } \forall (i, j) \in E& \label{k_transport_r}
\end{align}

{Solving the mathematical model formulated in this section using the exact approach is time-consuming. We therefore propose a tailored solution approach in Section \ref{Solution approach for preference-based multimodal trip planning} to reduce the computation time.}

\section{Solution approach for preference-based multi-modal trip planning}
\label{Solution approach for preference-based multimodal trip planning}

We first introduce the dynamic planning framework in Section \ref{Dynamic planning}, then we illustrate each type of transport mode and how they are synchronized in Sections \ref{Public Transport} to \ref{Private vehicles}.

\subsection{Dynamic planning}
\label{Dynamic planning}

The dynamic planning approach based on the rolling horizon is illustrated in Algorithm \ref{alg:Dynamic planning}. At each time step, the solutions considering all passengers are generated and recorded. The specific solution for each passenger is also recorded, including travel cost, time, utility, transport mode, route, and schedule. Algorithm inputs include vehicles \(K\) and initial locations of vehicles \(C\). Algorithm output is the (near) optimal solution \(Z\). For a new request at time \(t\), the active requests will be obtained, which include new requests as well as requests that have been scheduled but have not yet reached their intended destinations (line 4 in Algorithm \ref{alg:Dynamic planning}). Solutions \(Z_{\textrm{integrated}}^t\) that use integrated services with both PT and shared mobility, and solutions \(Z_{\textrm{micro}}^t\) that only use shared mobility are obtained using Algorithms \ref{alg:Integrated planning of PT and SM} and \ref{Micromobility_alg} (lines 5 and 6). Solutions \(Z_{\textrm{integrated}}^t\) and \(Z_{\textrm{micro}}^t\) are combined if these two solutions serve distinct requests (lines 7-9). Finally, The best solution \(Z^t\) is obtained according to Objective Function \eqref{obj2} (line 10). If there are unserved requests, the solution \(Z_{\textrm{private}}^t\) using private vehicles will be obtained by Algorithm \ref{Private} (lines 11-13). The solution set \(Z\) will be updated by adding \(Z^t\) or \(Z_{\textrm{private}}^t\) (line 14).

\begin{algorithm}[H]
\SetAlgoLined
 \textbf{Input:} $K$, $C$;
 \textbf{Output:} $Z$\tcp*{$Z$ represents the (near) optimal solution(s).}
 \While{True}{

 \If{{a new request $r$ is received at time $t$}}{
 obtain the active requests $R^t$\;
  $Z_{\textrm{integrated}}^t, C^t$ = $IntegratedPlanning(PT_{\textrm{schedule}}, K, C^{t-1}, R^t)$ (Algorithm \ref{alg:Integrated planning of PT and SM} in Section \ref{Public Transport})\;

 $Z_{\textrm{micro}}^t, C_{\textrm{micro}}^t$ = $Micromobility(K_{\textrm{micro}}, C_{\textrm{micro}}^{t-1}, R^t)$ (Algorithm \ref{Micromobility_alg} in Section \ref{Micromobility})\;
 \If{$Z_{\textrm{integrated}}^t$ or $Z_{\textrm{micro}}^t$ exists}{
 find all combinations $Z^t$ of $Z_{\textrm{micro}}^t$ and $Z_{\textrm{integrated}}^t$ that serve distinct requests in $R^t$
 }
 Find the solution(s) $Z^t$ with the maximum values of Objective Function \ref{obj2}\;
 \If{no solution can be found for unserved requests $R^t_{unserve}$}{
 $Z_{\textrm{private}}^t$ = $Private(R^t_{unserve}, K_{\textrm{private}})$ (Algorithm \ref{Private} in Section \ref{Private vehicles})
 }
 add $Z^t$ or $Z_{\textrm{private}}^t$ to solution set $Z$.
 }
  }
 \caption{Dynamic planning}\label{alg:Dynamic planning}
\end{algorithm}

The initial input of Algorithm \ref{alg:Dynamic planning} includes PT schedules $PT_{\textrm{schedule}}$ and SM information, such as transport mode, the type of vehicles (e.g., floating or dock-based micro-mobility), speed, capacity, cost, and initial locations of vehicles. When new requests are received, all available services are considered as potential alternatives for serving the request. The request can be served by the combination of PT and SM, or only one transport mode, such as sharing micro-mobility and private vehicles. When PT is considered, both access and egress trips are considered and each trip can be served by the SM. Transfers between SM are not considered. The specific approaches and illustrations for different types of services can be found in Sections \ref{Public Transport} to \ref{Private vehicles}.

\subsection{Integrated planning of PT and other transport modes}
\label{Public Transport}
Timetables of PT services are obtained from General Transit Feed Specification (GTFS) data \citep{GTFS}. GTFS data is a standardized format for PT data. It provides a structured and consistent way to describe transit schedules, routes, stops, and other relevant information for various modes of PT, such as buses, trains, and subways. In this study, we use routes, stops, schedules, trips, and fare information in GTFS data. {The PT schedules are predefined and the size of the PT fleet is fixed.} The capacity of PT is assumed to be unlimited. When both PT and SM services/private transport modes are considered, we combine their services to achieve seamless transport, as illustrated in Algorithm \ref{alg:Integrated planning of PT and SM}. The inputs include vehicles \(K\), requests \(R\), and latest locations of vehicles \(C^{t-1}\). The outputs include the solution with integrated services $Z_{\textrm{integrated}}^t$ and the new locations of vehicles $C^t$. For each request \(r\), first, the suitable PT services are obtained by identifying the nearby stations and matching the schedules with the request (lines 3-8). The request is divided by the PT service into three segments: the access request \(r_{\textrm{access}}\), PT request \(r_{\textrm{PT}}\), and egress request \(r_{\textrm{egress}}\) (line 9). The access and egress requests are added to the set of unplanned sub-requests \(R^{r}\) for request \(r\) (line 10). For each suitable PT service, there is a set of unplanned sub-requests \(R^{r}\). Assume that there are \(n_r\) PT services that can be used for request \(r\) and the platform aims to find the system optimal solution, then there will be \(\prod_{1}^{|R|} n_r\) combinations of sub-requests that need to be planned using shared mobility for all requests. We use \(R_{\textrm{subpool}}^{i} \in R_{\textrm{subpool}}\) to represent their combination. For each \(R_{\textrm{subpool}}^{i}\), first, we calculate the cost of PT; then, we calculate the cost for using each type of shared mobility and choose a solution \(Z^i\) according to Objective Function \ref{obj2} (lines 17-23); and finally, we choose the best \(Z_{\textrm{integrated}}^t\) from all \(Z^i\) for \(R_{\textrm{subpool}}^{i} \in R_{\textrm{subpool}}\) according to Objective Function \ref{obj2} (line 24).

\begin{algorithm}[h]
\SetAlgoLined
 \textbf{Input:} $PT_{\textrm{schedule}}$, $K$, $C^{t-1}$, $R^t$\;

 \textbf{Output:} $Z_{\textrm{integrated}}^t, C^t$\;

 obtain the current request $r$ from $R^t$; find PT stations $S_p$ and $S_d$ that are less than $n$ km from $p_r$ and $d_r$, respectively\;
 \For{$s_p \in S_p$}{
 \For{$s_d \in S_d$}{
 \For{${pt} \in K_{\textrm{PT}}$}{
 \If{${pt}$ operates between $s_p$ and $s_d$ and other transport modes can reach $s_p$/$d_r$ before/after ${pt}$ departs/arrives}{
 segment $r$ to access request $r_{\textrm{access}}$, PT request $r_{\textrm{PT}}$, and egress request $r_{\textrm{egress}}$ by departure and arrival stations/times of ${pt}$\;

 add $r_{\textrm{access}}$ and $r_{\textrm{PT}}$ to $R^{r}$\;
 }
 }
 }}

 obtain all possible request pools $R_{subpool}^{i} \in R_{subpool}$ by combining sub-requests of one ${pt}$ in $R^{r}$\;
\For{$R_{subpool}^{i} \in R_{subpool}$}{
calculate the cost of PT $F_{\textrm{PT}}^{i}$ for $R_{subpool}^{i}$\;
$Z_{\textrm{ride}}^i$, $C_{\textrm{ride}}^{t}$ = RidePooling($K_{\textrm{ride}}$, $C_{\textrm{ride}}^{t-1}$, $Z_{\textrm{ride}}^{t-1}$, $R_{subpool}^{i}$) (Algorithm \ref{RidePooling} in Section \ref{Ride pooling})\;
$Z_{\textrm{micro}}^i$, $C_{\textrm{micro}}^{t}$ = Micromobility($K_{\textrm{micro}}$, $C_{\textrm{micro}}^{t-1}$, $R_{subpool}^{i}$) (Algorithm \ref{Micromobility_alg} in Section \ref{Micromobility}) \;
$Z_{\textrm{private}}^i$ = Private($R_{subpool}^{i}$, $K_{\textrm{private}}$) (Algorithm \ref{Private} in Section \ref{Private vehicles})\;

Combine the PT and other transport modes to get solution $Z^i$ according to Objective Function \ref{obj2}.
}
choose the best solution $Z_{\textrm{integrated}}^t$ from all solutions $Z^i$ according to Objective Function \ref{obj2}.
 \caption{{Integrated Planning}}\label{alg:Integrated planning of PT and SM}
\end{algorithm}

\subsection{Routing algorithms for shared mobility}
\label{Shared mobility}
{Both ride-pooling and micro-mobility are considered SM services. This section introduces a tailored meta-heuristic algorithm for the routing of ride-pooling vehicles and micro-mobility vehicles.}

\subsubsection{{Adaptive Large Neighborhood Search for routing of a ride-pooling fleet}}
\label{Ride pooling}

Ride-pooling services can be offered by shared taxis, vans, or mini-buses with various capacities and speeds. Solving the ride pooling problem to optimality is computationally expensive \citep{santos2015taxi}, therefore we adopt Adaptive Large Neighborhood Search (ALNS) \citep{ropke2006adaptive} for ride pooling, as shown in Algorithm \ref{RidePooling}. According to the characteristics of the ride pooling problem, ALNS adapts its strategy to choose insertion and removal operators and achieve optimal or near-optimal solutions {by utilizing simulated annealing} (lines 4 and {16}) \citep{zhang2022synchromodal,ropke2006adaptive,zhang2022preference,azadeh2022choice}. The removal operator selectively removes passengers or ride requests from existing routes while maintaining feasibility (lines 5 and 8). It aims to optimize the ride pooling solution by identifying which passengers can be removed without compromising the quality of service, leading to more efficient vehicle routes and reduced travel times and costs. In contrast, the insertion operator adds new passengers or ride requests to existing vehicle routes (line 12). It evaluates potential insertion points within the current solution, considering constraints and optimization objectives. In insertion operators, new requests can be inserted in any position of the existing routes, and the visiting schedules can be adjusted as long as no scheduled passenger experiences a delay. ALNS will choose the solution with a better objective value between the current solution $Z_{\textrm{current}}$ and the last solution $Z_{\textrm{last}}$ (line 14). If the objective values of $Z_{\textrm{current}}$ and $Z_{\textrm{last}}$ are the same, the current solution will be assigned to $Z_{\textrm{last}}$ according to a probability generated by simulated annealing. Finally, the obtained solution in the current iteration $Z_{\textrm{last}}$ will be compared with the historical best solution $Z_{\textrm{ride}}^t$, and $Z_{\textrm{last}}$ will be assigned to $Z_{\textrm{ride}}^t$ if $Z_{\textrm{last}}$ is better than $Z_{\textrm{ride}}^t$ (line 15). The above steps are repeated until a predefined number of iterations is reached.

\begin{algorithm}[h]
\SetAlgoLined
 \textbf{Input:} $K_{\textrm{ride}}$, $C_{\textrm{ride}}^{t-1}$, $Z_{\textrm{ride}}^{t-1}$, $R_{subpool}^{i}$;
 \textbf{Output:} $Z_{\textrm{\textrm{ride}}}^t$\tcp*{$Z_{\textrm{\textrm{ride}}}^t$ represents the optimal solution(s) using ride pooling services.}
set $T_{\textrm{Temp}}>0$;

 $Z_{\textrm{last}} \leftarrow Z_{\textrm{ride}}^{t-1}$; $Z_{\textrm{\textrm{ride}}}^t \leftarrow Z_{\textrm{last}}$; $R_{\textrm{pool}} \leftarrow R_{subpool}^{i}$; \tcp*{$Z_{\textrm{last}}$ represents the solution in the last iteration.}

 \Repeat{the preset number of iterations has been completed}{
  at the beginning of each segment, update weights of operators and select operators depending on weights\;

$Z_{\textrm{current}} \leftarrow Z_{\textrm{last}}$; [$Z_{\textrm{current}}, R_{\textrm{pool}}$] = $RemovalOperator(Z_{\textrm{current}}, R_{\textrm{pool}})$; $flag = False$\;
  \While{requests in $R_{\textrm{pool}}$ have not been served}
  {
  \eIf{$flag==True$}{[$Z_{\textrm{current}}, R_{\textrm{pool}}$] = $RemovalOperator(Z_{\textrm{current}}, R_{\textrm{pool}})$}{$flag=True$}
  [$Z_{\textrm{current}}, R_{\textrm{pool}}$] = $InsertionOperator(Z_{\textrm{current}}, R_{\textrm{pool}})$\;

 }

Compare $Z_{\textrm{current}}$ and $Z_{\textrm{last}}$; update $Z_{\textrm{last}}$ by the simulated annealing \citep{ropke2006adaptive}.\;
$Z_{\textrm{\textrm{ride}}}^t \leftarrow Z_{\textrm{last}}$; $T_{\textrm{Temp}} \leftarrow T_{\textrm{Temp}} \cdot c $\tcp*{$c$ is the cooling rate.}
 }
 \caption{ALNS algorithm for ride-pooling vehicles}\label{RidePooling}
\end{algorithm}

ALNS explores different combinations of passenger-vehicle assignments, routes, and schedules by using a series of insertion and removal operators. {The selection of insertion and removal operators is critical for the effectiveness of the ALNS framework. For insertion operators, we utilize the Greedy Insertion, Random Insertion, Most Constrained First Insertion, and Regret Insertion operators. For removal operators, we employ the Worst Removal, Random Removal, Related Removal, and History Removal operators. Each operator is chosen based on its ability to address specific aspects of the optimization problem. Insertion operators such as the Greedy Insertion and Regret Insertion operators focus on minimizing immediate and future costs, respectively. Random Insertion and Most Constrained First Insertion operators ensure solution diversity and feasibility under tight constraints. Removal operators like Worst Removal and History Removal prioritize improving high-cost or suboptimal routes, while Random Removal and Related Removal operators introduce stochasticity and optimize clustered requests.} {The selected insertion and removal operators are described in detail below:}

\noindent\textbf{Greedy Insertion Operator}: This operator examines each feasible solution using a single vehicle or multiple vehicles and adds the request to the most optimal route(s) \citep{wolfinger2021large}. {This operator is based on the principle of minimizing immediate costs by selecting the most optimal route for a given request. It is particularly effective for improving initial solutions and finding feasible insertions with low computational effort. Its inclusion ensures that the algorithm can quickly converge to a high-quality solution in the early stages of optimization.}

\noindent\textbf{Random Insertion Operator}: It selects vehicles and insertion locations randomly and incorporates the request once a feasible solution is identified \citep{danloup2018comparison}. {This operator introduces randomness to diversify the search process, preventing premature convergence to local optima. It explores less obvious insertion points, helping the algorithm escape local optima and improve solution robustness.}

\noindent\textbf{Most Constrained First Insertion Operator}: Requests are prioritized based on a weighted evaluation of various factors such as the distance between origins and destinations, the loads, and the time windows. This operator systematically handles the most challenging requests first, potentially simplifying subsequent insertions \citep{danloup2018comparison}. {This operator handles the most challenging requests first, ensuring that the algorithm does not overlook critical constraints like tight time windows or high load demands. It is designed to tackle scenarios where certain requests might otherwise remain unassigned due to conflicts.}

\noindent\textbf{Regret Insertion Operator}: This method involves calculating a ``regret" value for not choosing the optimal route immediately by considering the difference in cost between the best current option and the next best. This foresight helps in minimizing future costs and inefficiencies by selecting routes that leave open the most beneficial future possibilities \citep{zhang2022synchromodal}. {This operator takes future costs into account by evaluating the regret value of not choosing the best route immediately, ensuring decisions made in the present do not lead to inefficiencies in subsequent iterations.}

\noindent\textbf{Worst Removal Operator}: It identifies and removes the request with the highest cost from routes \citep{wolfinger2021large}. {By removing the request with the highest cost, this operator prioritizes the improvement of high-cost routes. It ensures that resources are reallocated to minimize overall inefficiencies in the system.}

\noindent\textbf{Random Removal Operator}: This operator randomly removes a request from a subset of vehicles \citep{danloup2018comparison}. {Similar to the random insertion operator, this operator introduces stochasticity into the process. Removing requests randomly allows the algorithm to explore alternative configurations and potentially discover better solutions.}

\noindent\textbf{Related Removal Operator:} When a request is removed, it also removes similar requests, leveraging the shared attributes to optimize the removal process across multiple routes \citep{danloup2018comparison}. {This operator removes requests with shared attributes, such as overlapping time windows. It is particularly useful in identifying clusters of requests that might benefit from reorganization.}

\noindent\textbf{History Removal Operator}: It utilizes historical data to identify requests that may be sub-optimally placed. By focusing on requests whose current insertion costs are higher than their recorded lowest costs, this operator aims to re-optimize routes more cost-effectively \citep{zhang2022synchromodal}. {By leveraging historical cost data, this operator identifies suboptimal requests and repositions them to improve solution quality. It is valuable for iteratively refining solutions based on past performance.}

\subsubsection{{Heuristic algorithm for routing of micro-mobility fleet}}
\label{Micromobility}

Micro-mobility services are offered by means of scooters, bikes, or e-bikes, which can either be dock-based or floating. In the context of dock-based micro-mobility, vehicles are confined to movements between predefined stations, necessitating passengers to pick up and return vehicles at designated stations. Conversely, floating micro-mobility allows vehicles to be relocated to any location, with passengers having direct access to these vehicles for their journeys. {In contrast to ride pooling, where vehicles can travel to passengers' locations for pickups, micro-mobility vehicles are stationary and passengers are required to access micro-mobility vehicle locations before riding them.}

Algorithm \ref{Micromobility_alg} provides the pseudo-code for the micro-mobility services. The status of a micro-mobility vehicle is categorized as either ``idle" or ``busy". When a micro-mobility vehicle is chosen for a passenger's ride, it is promptly designated as ``busy" (line 19). At each time step, the status of the vehicle is subject to an update. If the vehicle has been released by the previous passenger, it is then reset to ``idle", signifying its availability for the next passenger's use (lines 4 to 6). For all requests $r \in R^t$, the approach iterates all available micromobility vehicles $k \in K_{\textrm{micro}}$. There are two types of micromobility services: dock-based and free-floating fleets. For a dock-based fleet, passengers must pick up and release vehicles at designated stations (line 9). For a free-floating fleet, passengers go directly to the vehicle and use it to reach their destination (line 11). Passengers first walk to the micro-mobility vehicle's location, board it for a ride, and then either park the vehicle before reaching their destination (in the case of a dock-based fleet) or ride it directly to their destination (for a free-floating fleet). We calculate the total travel time for option $k$, including both the riding and walking time. Suppose this total time is less than the longest allowable travel time $b_r - a_r$, the vehicle $k$ is considered as a feasible alternative (lines 13 to 14). We use the ``First come, First serve" strategy to provide micro-mobility services to requests in $R^t$ (line 18).

\vspace{-0.5cm}

\begin{algorithm}[h]
\SetAlgoLined
 \textbf{Input:} $K_{\textrm{micro}}$, $C_{\textrm{micro}}^{t-1}$, $R^t$;
 \textbf{Output:} $Z_{\textrm{micro}}^t$, $C_{\textrm{micro}}^{t}$ \tcp*{$Z_{\textrm{micro}}$ represents the optimal solution using micro-mobility.}

  obtain the current request $r$ from $R^t$\;
 \For{$k \in K_{\textrm{micro}}$}{
 \If{status of $k$ is not ``idle" and $t > t^k$}{
 set status of $k$ as ``idle"\;
 }
 \If{status of $k$ is ``idle"}{
 \uIf{type of $k$ if dock-based}{
 find the closest station of $k$ from $d_r$ and set it as ${l'}^k$\;
 }
 \Else{
 ${l'}^k \leftarrow d_r$\;
 }
 \If{(($d_{p_r l^k}^k$ + $d_{{l'}^k d_r}^k$)/$v_r$ + $d_{l^k {l'}^k}^k$/$v_k$) $<$ $b_r - a_r$}{
 calculate utility $V_r$ using Equation \eqref{utility_from_system} and add $V_r$ and $k$ to the set of utilities $\textbf{V}_r$ for $r$\;
 }
 }
 }
 choose $k$ with the maximum utility in $\textbf{V}_r$\ and add it to the solution $Z_{\textrm{micro}}^t$\;
 update vehicle $K$'s location in $C_{\textrm{micro}}^{t}$ and change $k$'s status as ``busy"\;

 \caption{Micro-mobility vehicles}\label{Micromobility_alg}
\end{algorithm}

\vspace{-0.5cm}

\subsection{Privately owned vehicle}
\label{Private vehicles}
For private vehicles, private bikes, private scooters, and private cars are considered, as illustrated in Algorithm \ref{Private}. Walking is also considered as an alternative. The walking speed depends on the age \citep{walkingspeed,alves2020walkability}, while the average speeds for private vehicles are obtained from the literature and reports \citep{cyclingspeed}. If the traveling time $d_{p_r, d_r}^k / v_k$ ($v_k$ is speed) using private vehicle $k$ is shorter than the longest allowable travel time $b_r - a_r$, $k$ will be considered as an alternative. The alternative with maximum utility will be chosen and added to the solution using private vehicles $Z_{\textrm{private}}^t$.

\begin{algorithm}[h]
\SetAlgoLined
 \textbf{Input:} $K_{\textrm{private}}$, $R^t$;
 \textbf{Output:} $Z_{\textrm{private}}^t$ \tcp*{$Z_{\textrm{private}}^t$ represents the optimal solution(s) using private vehicles.}
  obtain the current request $r$ from $R^t$\;

 \For{$k \in K_{\textrm{private}}^r$}{
 \If{$d_{p_r, d_r}^k / v_k< b_r - a_r$}{
  calculate utility $V_r$ using Equation \ref{utility_from_system} and add $V_r$ and $k$ to the set of utilities $\textbf{V}_r$ for $r$\;
 }
 }
 choose $k$ with the maximum utility in $\textbf{V}_r$\ and add it to the solution $Z_{\textrm{private}}^t$\;

 \caption{Privately owned vehicles}\label{Private}
\end{algorithm}

\section{{Case study}}
\label{Case Study}
In the following, we present our experiment results. Initially, the setup of the experiments is outlined in Section \ref{Experimental setup}. Next, Section \ref{Solution method performance} presents a comparative analysis of the results obtained from the MIP model and the proposed meta-heuristic algorithm. Finally, in Section \ref{Numerical experiments and managerial insights}, we describe the numerical experiments conducted and {discuss} the managerial insights derived from these findings.

\subsection{{Experimental setup}}
\label{Experimental setup}

{Our research is conducted as part of the European Horizon project, Seamless Shared Urban Mobility (SUM), which provides the demand data and realistic transport network to this study. In this project, we collaborate with a real transport planner in the Netherlands, RET \citep{RET}, which is developing platforms to provide personalized multimodal trip plans to individual users. To simulate the system, we generate trip requests based on aggregated data obtained from surveys \citep{montes2023shared}. This data includes detailed choice components that are carefully tailored to the case study presented in this work.}

{We apply our model to a case study in the suburban area of Maassluis (Rotterdam), the Netherlands. Compared to the Rotterdam city center, Maassluis has lower population densities, greater distances between destinations, and less extensive PT infrastructure. As a result, residents in the Maassluis area may rely more heavily on privately owned cars for their transport needs, leading to increased traffic congestion, environmental impact, and commuting stress. This presents an opportunity to explore innovative multi-modal solutions that integrate existing PT services with SM services. According to a commonly used Hexagonal Hierarchical Geospatial Indexing System \citep{H3}, the studied area is divided into hexagonal grid cells.}

\begin{figure}[H]
\centering
\begin{subfigure}[b]{0.495\textwidth}
\includegraphics[width = 0.8\textwidth, height=5.5cm]{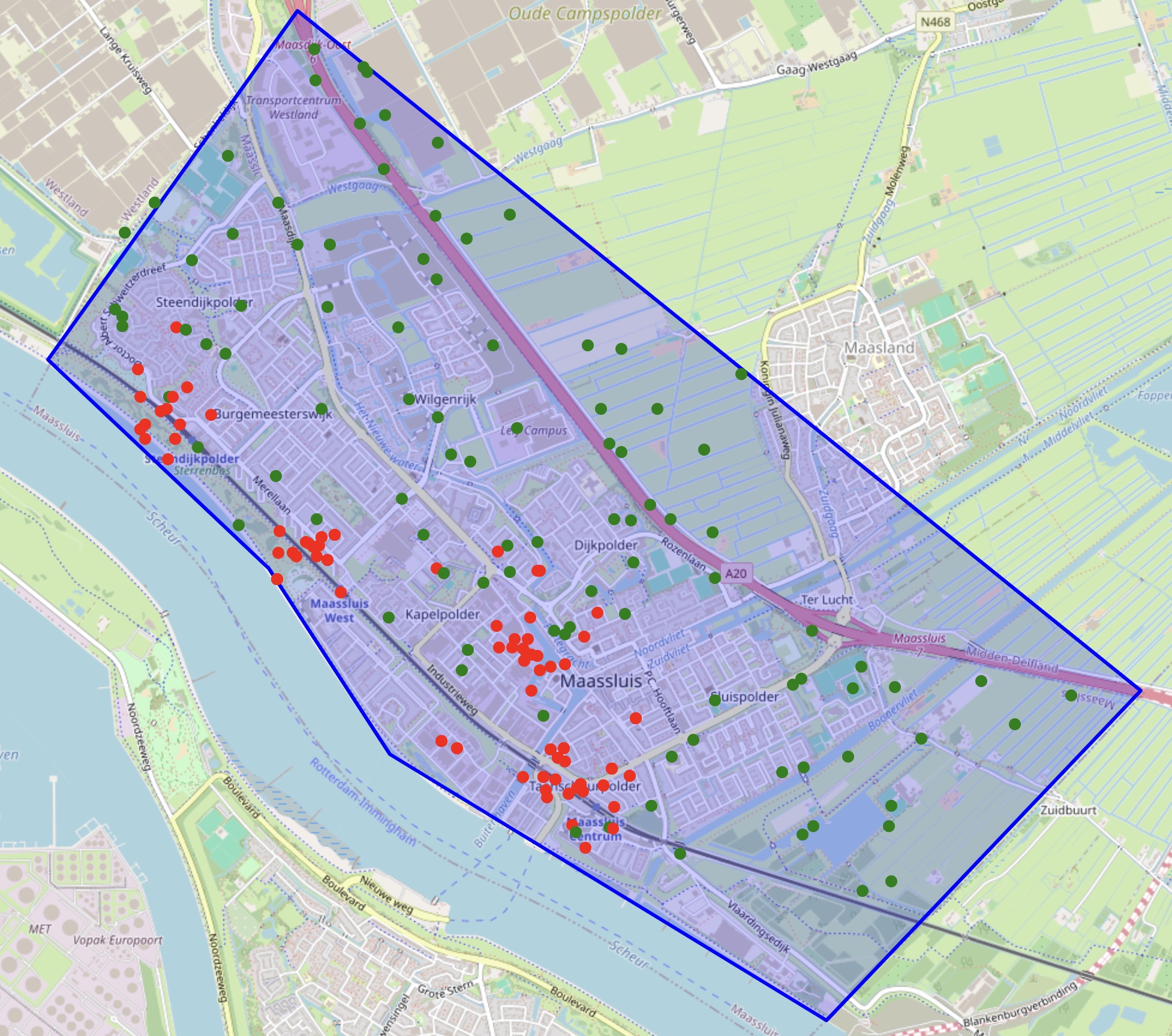}
    \caption{Origins (green color) and\\ destinations (red color) of passengers.}
    \label{fig:Request_on_map.}
\end{subfigure}
\begin{subfigure}[b]{0.495\textwidth}
    \includegraphics[width = 0.8\textwidth, height=5.5cm]{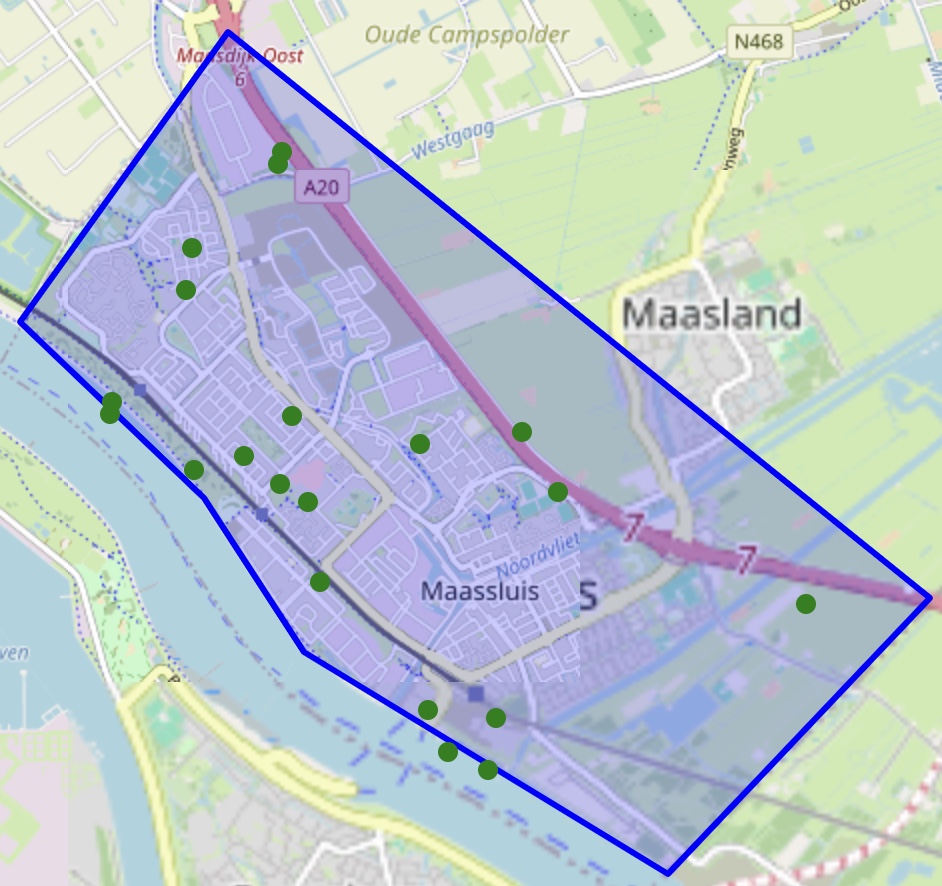}
    \caption{Initial locations of scooters.\\\text{ }}
    \label{fig:Scooter_on_map.}
\end{subfigure}

\begin{subfigure}[b]{0.495\textwidth}
 \includegraphics[width = 0.8\textwidth, height=5.5cm]{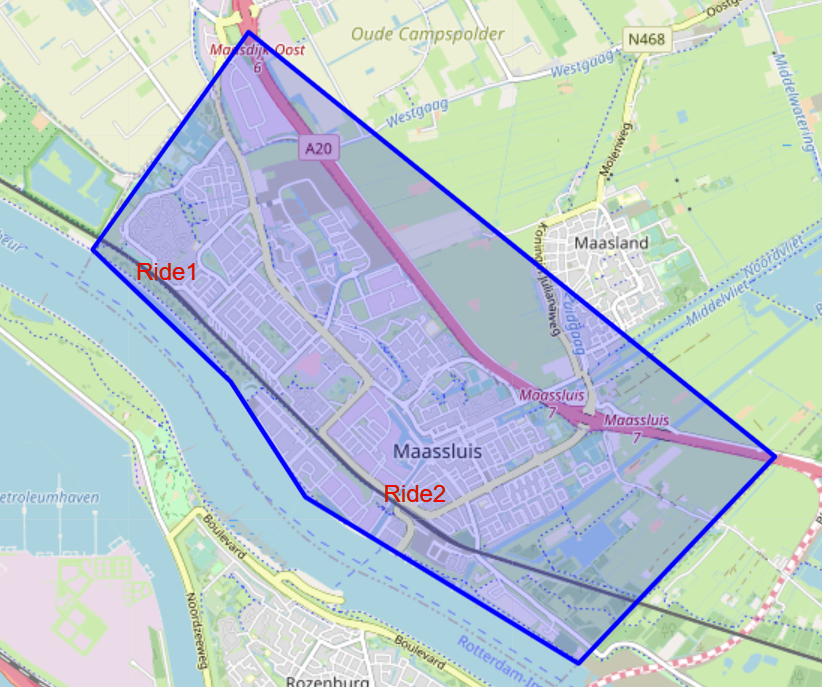}
    \caption{Initial locations of ride-pooling vehicles.}
    \label{fig:Ride_on_map.}
\end{subfigure}
\begin{subfigure}[b]{0.495\textwidth}
\includegraphics[width = 0.8\textwidth, height=5.5cm]{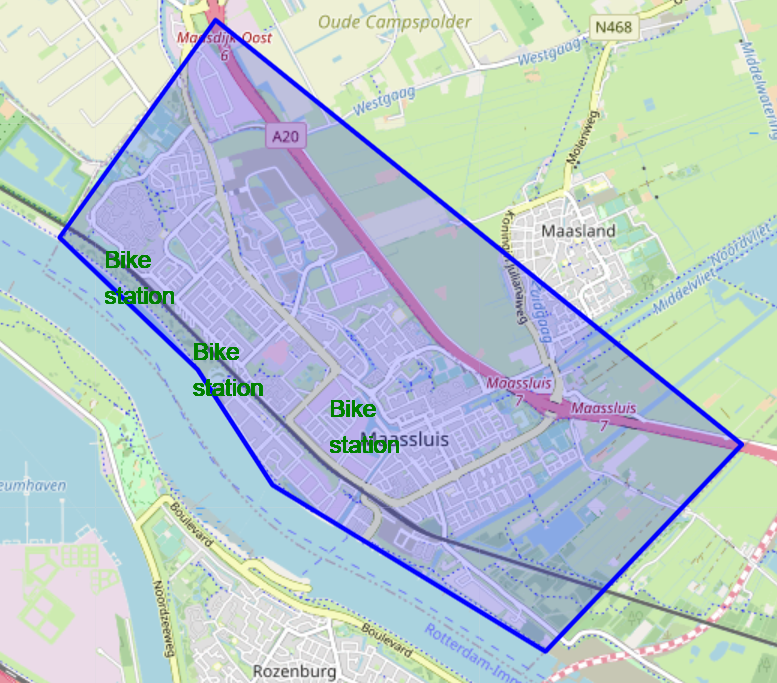}
    \caption{Bike stations.}
    \label{fig:Bike_on_map.}
\end{subfigure}
    \caption{Requests and vehicles in Maassluis area.}
    \label{fig:Requests and vehicles in Maassluis area.}

\end{figure}

As shown in Figure \ref{fig:Request_on_map.}, the origins of passengers are randomly generated in grid cells, and the destinations of passengers are generated based on a Gamma distribution, as outlined in the illustrations referenced in \cite{soza2024shareability}. {The Gamma distribution is defined by two parameters: the shape parameter ($k_c$) and the scale parameter ($s_c$). The shape parameter controls the clustering density of destinations around PT stations, with lower values indicating more concentrated clustering and higher values representing greater dispersion. The scale parameter adjusts the overall spread of destinations across the study area. This approach ensures realistic modeling of destination clustering around PT stations, reflecting typical spatial travel patterns.} The parameters for the utility functions are listed in Tables \ref{table:utility parameters} and \ref{table:utility parameters_beta} in Appendix \ref{Parameters in the utility function}. {The utility functions incorporate parameters for travel time, cost, and user-specific attributes, such as sensitivity to waiting and in-vehicle time. These parameters are based on stated preference surveys conducted in Rotterdam \citep{montes2023shared}. For example, the sensitivity to travel cost is higher for shared micro-mobility compared to other modes, indicating that users of shared micro-mobility are more price-sensitive. Similarly, waiting time ($\beta_{\text{PTwait}}$) has a stronger disutility for public transport users, reflecting their lower tolerance for delays. These parameters are essential for realistically modeling mode choice behavior and capturing the trade-offs users make between different attributes.} The time scale of experiments is a whole day. When designing instances, the time interval between passengers follows time-dependent Poisson distributions. During peak hours and non-peak hours, the $\lambda_t$ parameter in the Poisson distribution is set to 10 and 5, respectively. {The maximum expected travel time for a passenger} is randomly generated between 30 minutes to 120 minutes.

The distances between locations are obtained from OpenStreetMap \citep{OSM}. PT service schedules are obtained from GTFS data \citep{GTFS}. In total, we consider 367 PT scheduled trips, each comprising multiple stops. The capacity of PT is assumed to be unlimited. The fare information is sourced from Rotterdam's primary PT service providers, namely NS \citep{NS} and RET \citep{RET}. Figures \ref{fig:Scooter_on_map.}-\ref{fig:Bike_on_map.} show the locations of shared mobility. The stations of shared bikes are located in high-population-density areas. The initial locations of scooters are generated by drawing from a uniform distribution within the boundaries of the designated area. For the ride-pooling vehicles, their depots are strategically located close to PT stations. Various modes of transport have distinct speeds \citep{walkingspeed,bikespeed,scooterspeed}.

\subsection{{Evaluation of the solution method's performance}}
\label{Solution method performance}

{All experiments are conducted using Python 3.9 and executed on a Linux system equipped with 62 GB of RAM and an Intel Xeon E5 processor. {As the objective is to maximize social welfare, the proposed approach ensures that services are selected to achieve this goal. To validate the effectiveness of the approach, we compare the proposed heuristic algorithm (HA) with an exact approach (MIP solved by Gurobi). Table \ref{table:exact} shows the comparison between the MIP and HA in terms of quality of solutions and computation time.} As the {sizes of SM fleet and requests} scale up, the computation time of the exact approach experiences a notable increase. In contrast, the proposed heuristic algorithms efficiently yield equivalent solutions with significantly reduced computation times by more than 90\% in all instances. The heuristic algorithm can find the optimal solution in all cases, and it also finds better solutions for the cases where the MIP does not find the optimal solution within the allowable time (three hours). {When multiple services are available, the proposed approach selects the one that yields higher social welfare.} For the instance with 5 shared bikes, 5 shared scooters, and 2 ride-pooing vehicles, the exact approach cannot find a feasible solution within the time limit, while the heuristic algorithm can find a near-optimal solution in less than three minutes.}

\begin{table}[h]
\centering
\begin{adjustbox}{max width=\textwidth}
\begin{threeparttable}
\caption{The comparison between the MIP and heuristic algorithm}
\label{table:exact}
\begin{tabular}{c c c c c c c c c c}
 \toprule
 &&&&\multicolumn{2}{c}{$R_{\textrm{serve}}$} &\multicolumn{2}{c}{Social Welfare}&\multicolumn{2}{c}{Avg. CPU (s)}\\
 \cmidrule(lr){5-6} \cmidrule(lr){7-8}\cmidrule(lr){9-10}
  $K_{\textrm{bike}}$&$K_{\textrm{scooter}}$&$K_{\textrm{ride}}$&$R$&MIP& HA&MIP& HA& MIP& HA \\
\midrule
1&1&0&1&1&1&-3.8&-3.8&10&1\\

1&1&0&3&3&3&-13.3&-13.3&89&9 \\

1&1&0&5&5&5&-19.6&-19.6&549&25 \\

1&1&0&7&7&7&-29.3&-29.3&2485&39 \\

1&1&0&9&9&9&-33.9&-33.9&3704&42\\
\hline

1&5&0&1&1&1&-3.4&-3.4&98&2\\

1&5&0&3&3&3&-8.5&-8.5&634&7 \\

1&5&0&5&5&5&-16.6&-16.6&1353&15 \\

1&5&0&7&7&7&-23.0&-23.0&5485&36\\

1&5&0&9&9&9&-27.9&-27.9&7600&46 \\
\hline

5&1&0&1&1&1&-3.8&-3.8&79&4\\

5&1&0&3&3&3&-12.1&-12.1&609&26\\

5&1&0&5&5&5&-18.3&-18.3&1949&74 \\

5&1&0&7&7&7&-28.1&-28.1&6435&94 \\

5&1&0&9&\textbf{8}&\textbf{9}&\textbf{-39.5}&\textbf{-39.3}&10800*&146 \\
\hline
5&5&0&1&1&1&-3.4&-3.4&430&3 \\

5&5&0&3&3&3&-8.5&-8.5&3109&29\\

5&5&0&5&5&5&-16.6&-16.6&6949&91\\

5&5&0&7&7&7&-23.0&-23.0&9435&113 \\

5&5&0&9&9&9&\textbf{-28.2}&\textbf{-27.9}&10800*&175 \\
\hline
1&1&2&1&1&1&-3.8&-3.8&90&2 \\

1&1&2&3&3&3&-13.4&-13.4&489&4 \\

1&1&2&5&5&5&-19.4&-19.4&1549&8 \\

1&1&2&7&7&7&-29.0&-29.0&6485&40 \\

1&1&2&9&9&9&-37.9&-37.9&8600&73 \\
\hline
5&5&2&1&1&1&-3.8&-3.4&930&5 \\
5&5&2&3&3&3&-12.1&-8.5&5009&27 \\

5&5&2&5&5&5&-21.3&-16.6&9949&83\\

5&5&2&7&\textbf{6}&\textbf{7}&\textbf{-31.6}&\textbf{-20.4}&10800*&103 \\

5&5&2&9&\textbf{-}&\textbf{9}&\textbf{-}&\textbf{-37.3}&10800*&157 \\

\bottomrule

\end{tabular}
\begin{tablenotes}
      \footnotesize
      \item $K_{\textrm{bike}}$, $K_{\textrm{scooter}}$, and $K_{\textrm{ride}}$ denote the numbers of shared bikes, shared scooters, and ride-pooling vehicles, respectively; $R$ and $R_{serve}$ denote the numbers of received requests and served requests, respectively;
      \item $\ast$ time limit reached (three hours). $-$ no solution has been found due to the time limitation.

    \end{tablenotes}
  \end{threeparttable}
  \end{adjustbox}
\end{table}

\begin{figure}[h]
    \centering
    \includegraphics[width=1\linewidth]{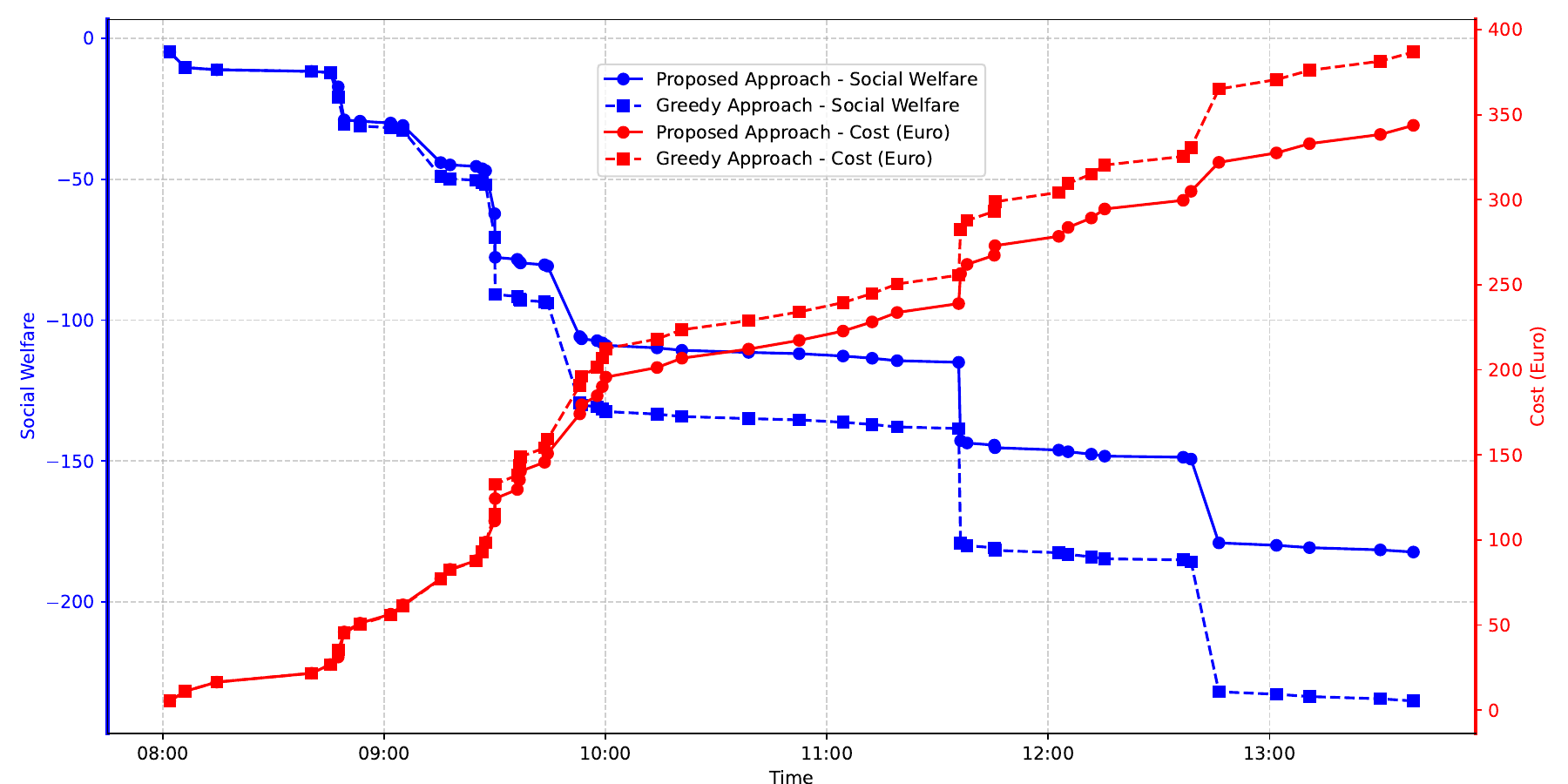}
    \caption{{Comparison of the proposed heuristic algorithm and a greedy approach.}}
    \label{fig:Comparison of the proposed heuristic algorithm and a greedy approach}
\end{figure}

{We also compare the performance of the proposed heuristic algorithm with a greedy approach on an instance involving 367 PT services, 20 SM services, and 50 requests over a time period from 8:00 AM to 2:00 PM. The greedy approach searches the solution space and selects possible alternatives in a locally optimal manner using greedy insertion, without exploring alternative placements through random insertion or considering regret values for different alternatives. We evaluate the performance of both approaches in terms of total social welfare and cost. In Figure \ref{fig:Comparison of the proposed heuristic algorithm and a greedy approach}, each marker on the line represents a single request. The results indicate that both approaches yield the same cost and social welfare for the first five requests. However, as the number of requests increases, the proposed approach consistently outperforms the greedy approach, with the gap widening over time. After serving the final request, the proposed approach achieves a 22\% improvement in social welfare and an 11\% reduction in cost compared to the greedy approach. The trend suggests that as demand grows, the benefits of the proposed approach become more pronounced, reinforcing the importance of intelligent trip planning. The proposed approach further illustrates that a well-designed optimization algorithm can effectively reduce costs while enhancing social welfare. To maximize social welfare, the multimodal transport platform should adopt a strategy similar to the proposed approach, which incorporates better decision-making mechanisms.}

\subsection{Numerical experiments and managerial insights}
\label{Numerical experiments and managerial insights}

{The following sub-sections present a detailed analysis of the experimental results {of large instances}. Section \ref{Modal shares of heterogeneous passengers} analyzes the modal shares of heterogeneous passengers. Sections \ref{Analysis of factors influencing multimodal trips}, \ref{Analysis of factors influencing single-mode trips}, and \ref{Analysis of factors influencing travel time and waiting time} specifically analyze the impact of factors that influence multi-modal trips, single-mode trips, and journey time, respectively.} The experiments are designed to {evaluate} the impacts and implications of different factors in multi-modal transport with PT and SM, including:

\noindent \textbf{User heterogeneity}: {The following segments of passengers are considered: A. the entire population; B1. commuters and B2. non-commuters; C1. passengers who are familiar with SM, and C2. passengers who are unfamiliar with SM; D1. passengers who use PT once per week or less, and D2. passengers who use PT more than once per week; E1. passengers with a low income, and E2. passengers with a high income; F1. male passengers, and F2. female passengers; G1. passengers with a low level of education, and G2. passengers with a high level of education; H1. passengers who have used SM, and H2. passengers who have not used SM; I1. passengers under 35 years old, I2. passengers between 35 and 65 years old, and I3. passengers over 65 years old.} {The parameters used for the utility function of different segments of passengers can be found in Appendix \ref{Parameters in the utility function}.}

    \noindent \textbf{Supply-demand balance}: We examine the performance of the proposed model under varying supply and demand conditions. These scenarios encompass a wide range of supply-demand ratios, from situations with abundant supply to cases of scarcity. {To reflect Rotterdam's suburban region's practical capacity limits}, the numbers of shared scooters and bikes are listed as [0, 5, 10, 50], while the number of requests varies within [5, 10, 30, 50, 100]. Given our case study's focus on a suburban area, ride-pooling vehicles are kept at low quantities, specifically within the range [0, 1, 2].

\noindent \textbf{Cost sensitivity}: {By adjusting the costs of micromobility services within the range of [1, 4, 8] (euros), we investigate how these varying costs influence mode choices and the overall efficiency of the system.}

\subsubsection{{Modal shares of heterogeneous passengers}}
\label{Modal shares of heterogeneous passengers}

\begin{figure}[H]
\centering
\begin{subfigure}[b]{0.283\textwidth}
\includegraphics[width = 0.8\textwidth]{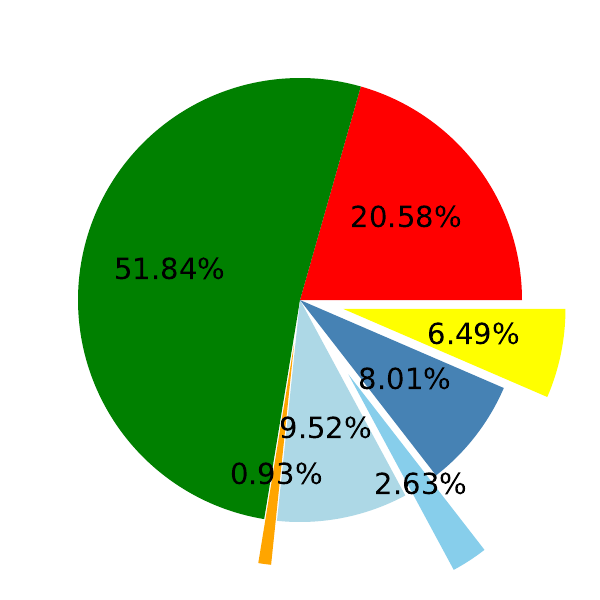}
    \caption{A. aggregated}
    \label{fig:Heuri1000end55000modal_sharesgeneral.}
\end{subfigure}
\hspace{-10mm}
\begin{subfigure}[b]{0.283\textwidth}
    \includegraphics[width = 0.8\textwidth]{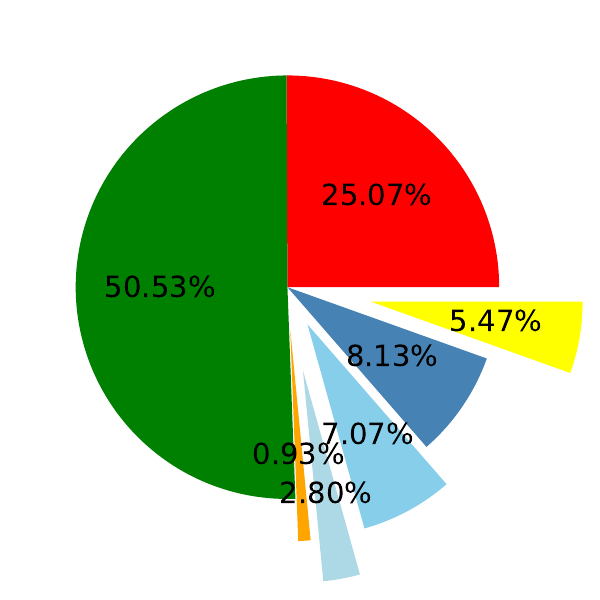}
    \caption{B1. commute}
    \label{fig:Heuri1000end55000modal_sharescommute}
\end{subfigure}
\hspace{-10mm}
\begin{subfigure}[b]{0.283\textwidth}
\includegraphics[width = 0.8\textwidth]{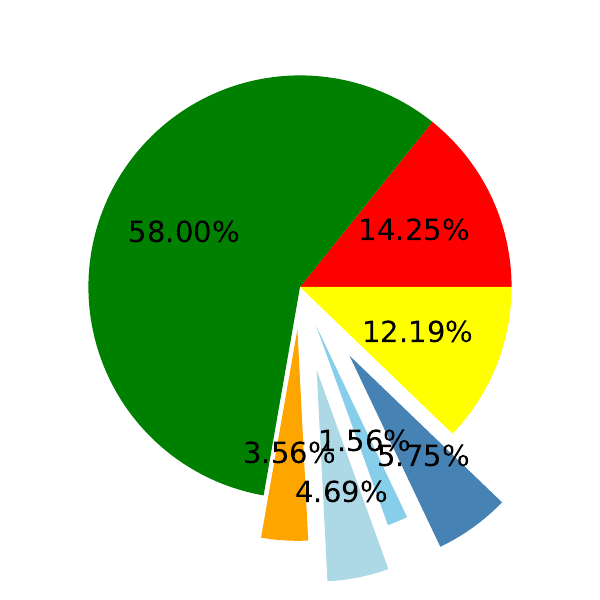}
    \caption{B2. not commute}
    \label{fig:Heuri1000end55000modal_sharesnot commute}
\end{subfigure}
\hspace{-10mm}
\begin{subfigure}[b]{0.283\textwidth}
\includegraphics[width = 0.8\textwidth]{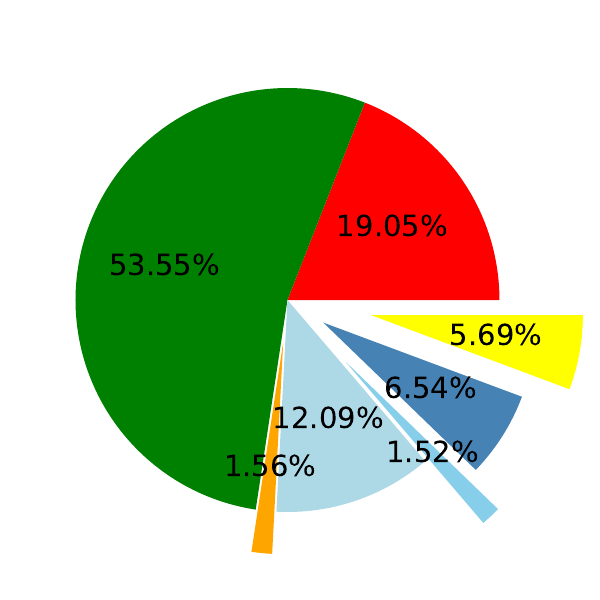}
    \caption{C1. familiar with SM}
    \label{fig:Heuri1000end55000modal_sharesfamiliarwith SM.}
\end{subfigure}
\hspace{-10mm}
\begin{subfigure}[b]{0.283\textwidth}
 \includegraphics[width = 0.8\textwidth]{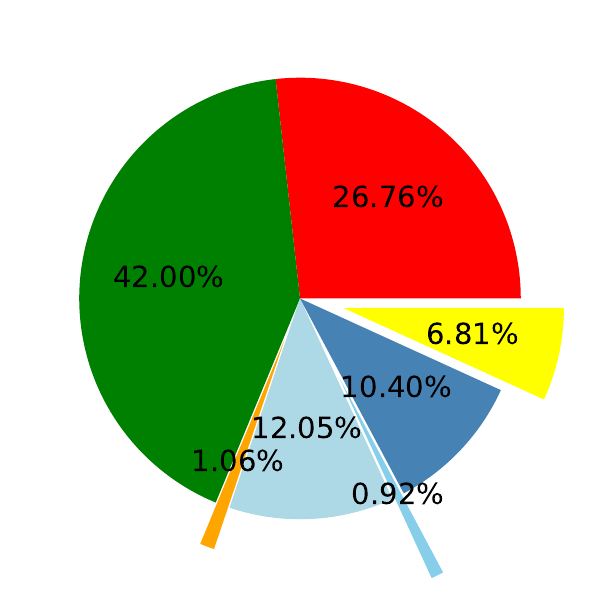}
    \caption{C2. not \\familiar with SM}
    \label{fig:Heuri1000end55000modal_sharesnot familiarwith SM}
\end{subfigure}
\hspace{-10mm}
\begin{subfigure}[b]{0.283\textwidth}
\includegraphics[width = 0.8\textwidth]{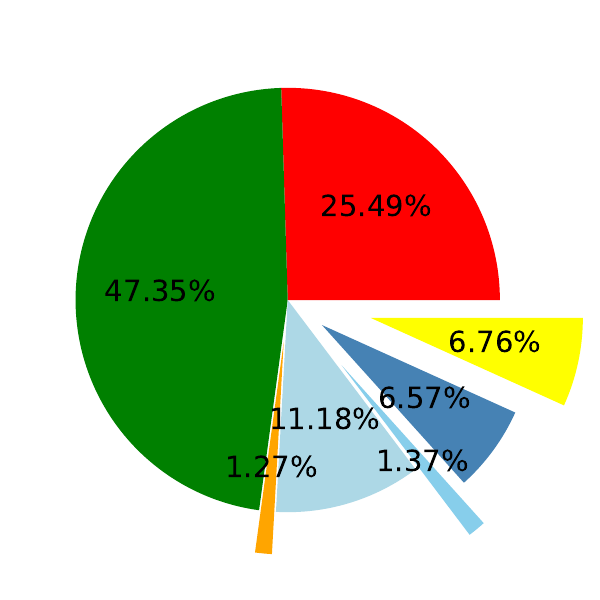}
    \caption{D1. use PT \\$\leq$ once/week}
    \label{fig:Heuri1000end55000modal_sharesuse PT  small onceweek.}
\end{subfigure}
\hspace{-10mm}
\begin{subfigure}[b]{0.283\textwidth}
\includegraphics[width = 0.8\textwidth]{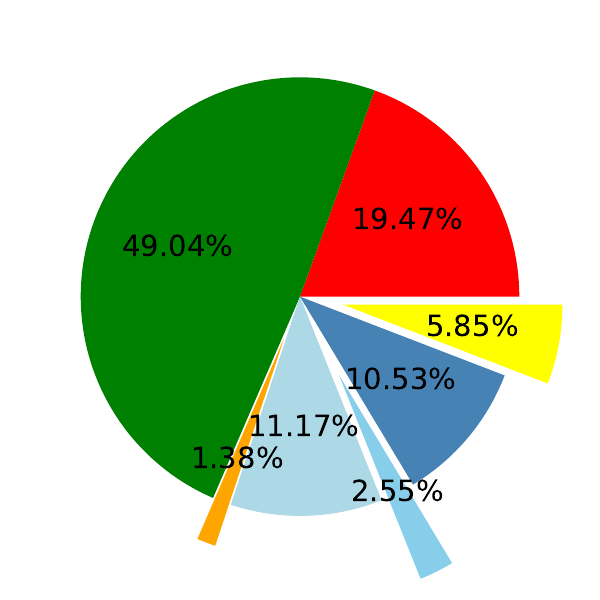}
    \caption{D2. use PT \\$>$ once/week}
    \label{fig:Heuri1000end55000modal_sharesuse PTlargeonceweek.}
\end{subfigure}
\hspace{-10mm}
\begin{subfigure}[b]{0.283\textwidth}
\includegraphics[width = 0.8\textwidth]{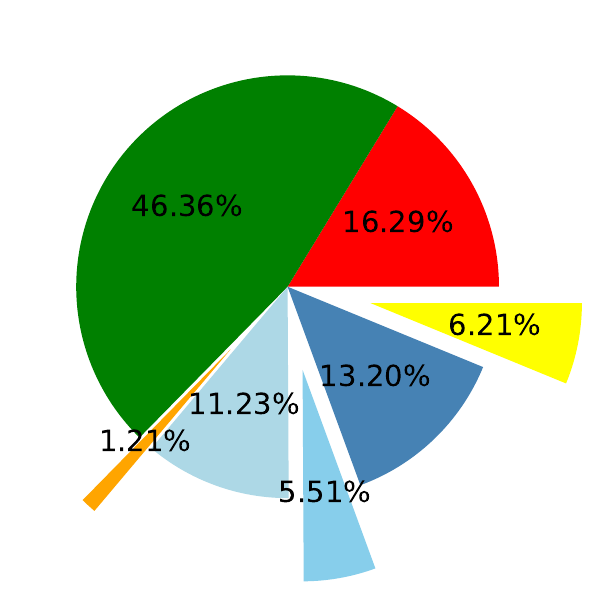}
    \caption{E1. low \\income}
    \label{fig:Heuri1000end55000modal_shareslowincome.}
\end{subfigure}
\hspace{-10mm}
\begin{subfigure}[b]{0.283\textwidth}
\includegraphics[width = 0.8\textwidth]{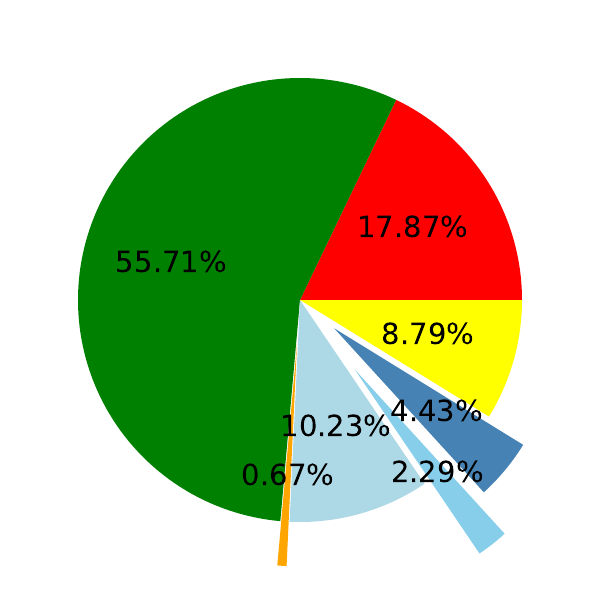}
    \caption{E2. high income}
    \label{fig:Heuri1000end55000modal_shareshighincome.}
\end{subfigure}
\hspace{-10mm}
\begin{subfigure}[b]{0.283\textwidth}
\includegraphics[width = 0.8\textwidth]{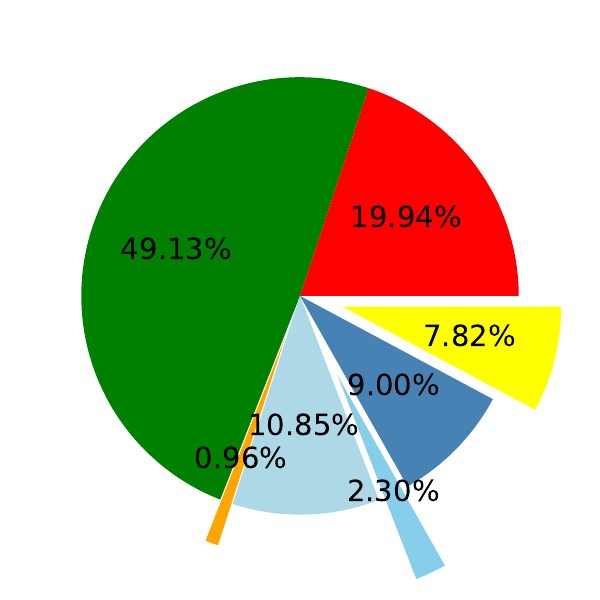}
    \caption{F1. male}
    \label{fig:Heuri1000end55000modal_sharesmale.}
\end{subfigure}
\hspace{-10mm}
\begin{subfigure}[b]{0.283\textwidth}
\includegraphics[width = 0.8\textwidth]{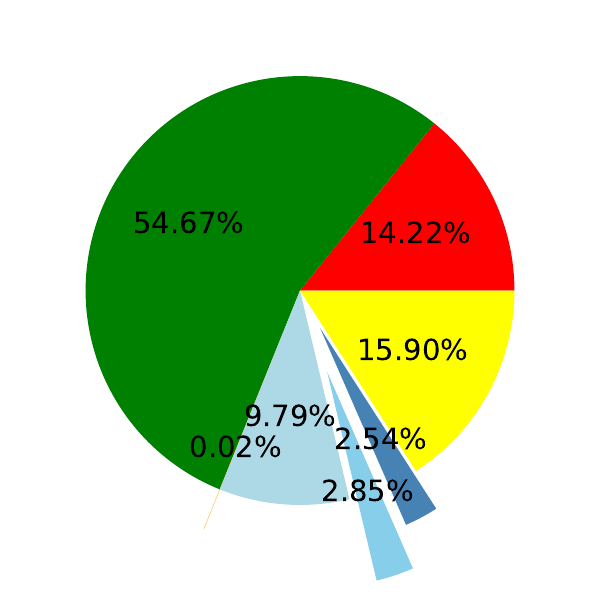}
    \caption{F2. female}
    \label{fig:Heuri1000end55000modal_sharesfemale.}
\end{subfigure}
\hspace{-10mm}
\begin{subfigure}[b]{0.283\textwidth}
\includegraphics[width = 0.8\textwidth]{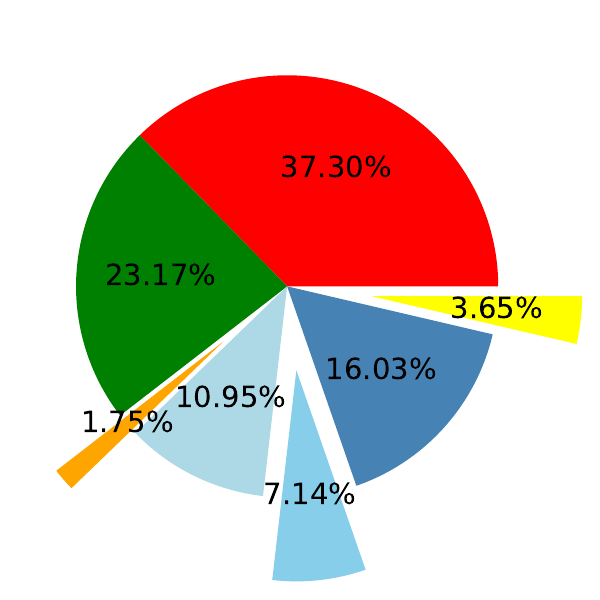}
    \caption{G1. low education}
    \label{fig:Heuri1000end55000modal_sharesloweducation.}
\end{subfigure}
\hspace{-10mm}
\begin{subfigure}[b]{0.283\textwidth}
\includegraphics[width = 0.8\textwidth]{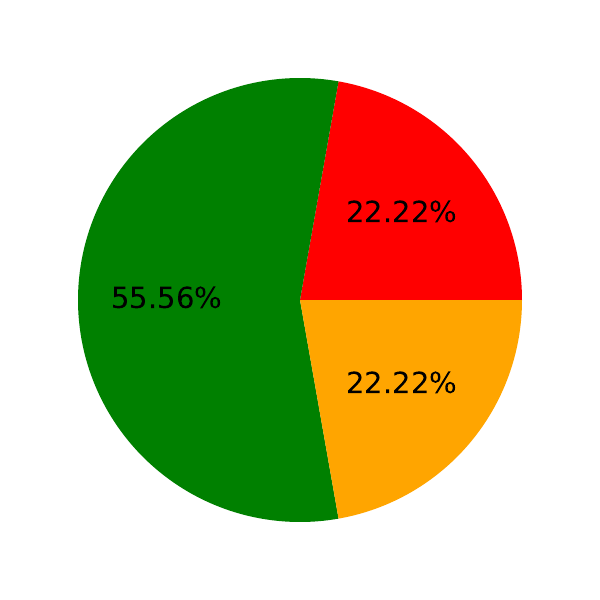}
    \caption{G2. high education}
    \label{fig:Heuri1000end55000modal_shareshigheducation.}
\end{subfigure}
\hspace{-10mm}
\begin{subfigure}[b]{0.283\textwidth}
\includegraphics[width = 0.8\textwidth]{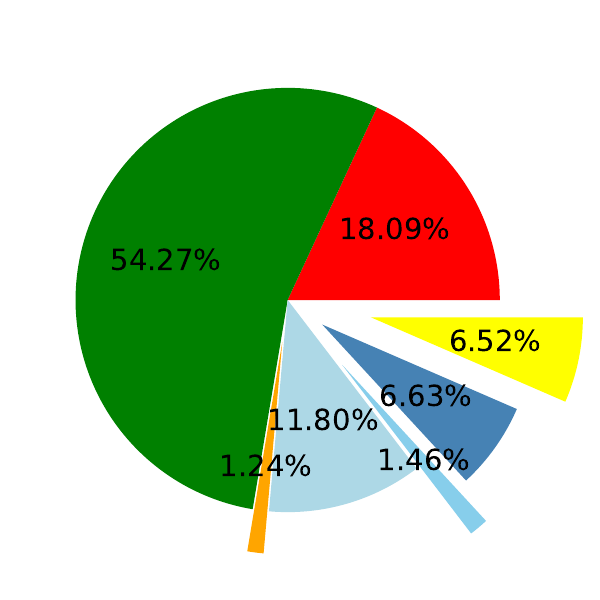}
    \caption{H1. used SM}
    \label{fig:Heuri1000end55000modal_shareshaving used SM.}
\end{subfigure}
\hspace{-10mm}
\begin{subfigure}[b]{0.283\textwidth}
\includegraphics[width = 0.8\textwidth]{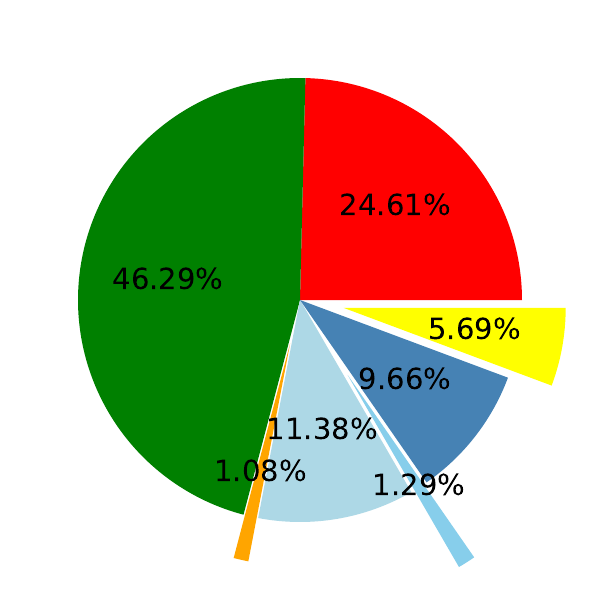}
    \caption{H2. have not used SM}
    \label{fig:Heuri1000end55000modal_shareshaving notused SM.}
\end{subfigure}
\hspace{-10mm}
\begin{subfigure}[b]{0.283\textwidth}
\includegraphics[width = 0.8\textwidth]{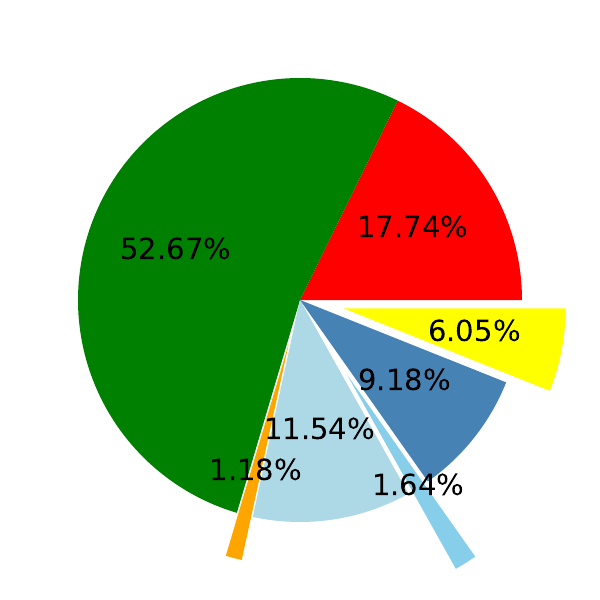}
    \caption{I1. age $\leq$ 35}
    \label{fig:Heuri1000end55000modal_sharesage small 35.}
\end{subfigure}
\hspace{-10mm}
\begin{subfigure}[b]{0.283\textwidth}
\includegraphics[width = 0.8\textwidth]{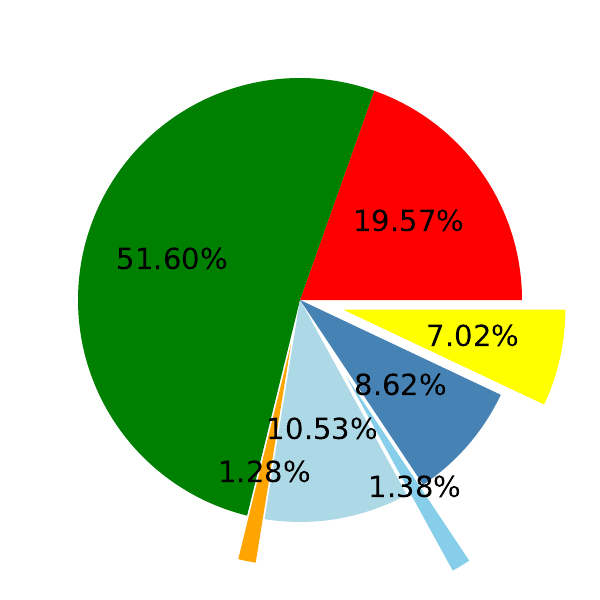}
    \caption{I2. 35 $<$ age $<$ 65}
    \label{fig:Heuri1000end55000modal_shares35  small age  small 65.}
\end{subfigure}
\hspace{-10mm}
\begin{subfigure}[b]{0.283\textwidth}
\includegraphics[width = 0.8\textwidth]{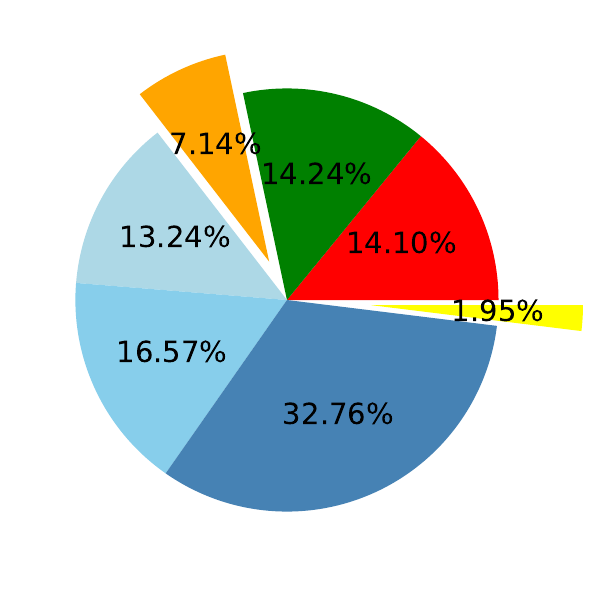}
    \caption{I3. age $\geq$ 65}
    \label{fig:Heuri1000end55000modal_sharesagelarge65.}
\end{subfigure}
\hspace{-10mm}
\begin{subfigure}[b]{0.283\textwidth}
\includegraphics[width = 0.8\textwidth]{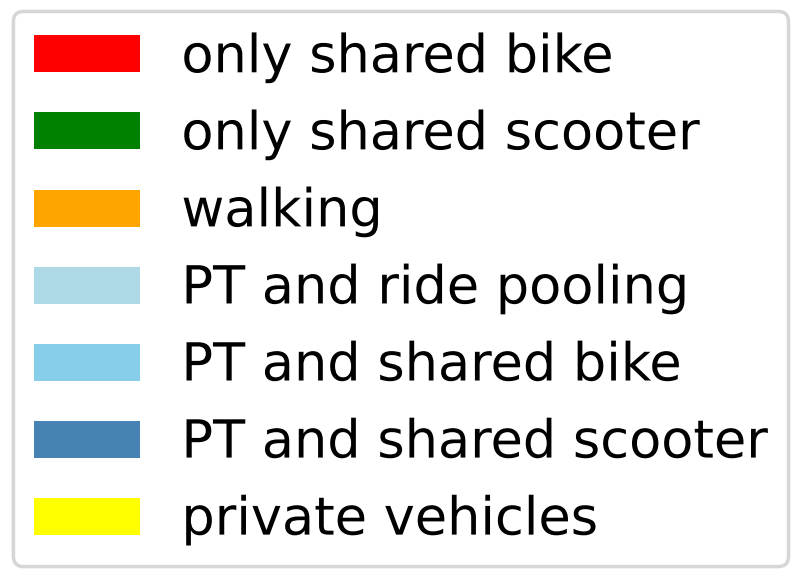}
    \caption*{}
    \label{fig:modal_share_legend.pdf}
\end{subfigure}
    \caption{Modal shares of different segments of passengers.}
    \label{fig:Modal shares of different segments of passengers.}
\end{figure}

Figure \ref{fig:Modal shares of different segments of passengers.} illustrates the distribution of modal shares among various passenger segments. Figure \ref{fig:Heuri1000end55000modal_sharesgeneral.} provides an overview of general passengers. Scooters dominate the modal share at 52\%, followed by bikes at 21\%, with walking representing the smallest share. Comparing commuters to non-commuters across Figures \ref{fig:Heuri1000end55000modal_sharescommute} to \ref{fig:Heuri1000end55000modal_sharesnot commute}, commuters exhibit an 11\% higher modal share for bikes (all differences mentioned are in percentage-points) and an 8\% higher share of PT and micromobility combinations. Conversely, commuters show a lower usage of PT and ride pooling, walking, and scooters by 2\%, 3\%, and 7\%, respectively. {Therefore, surplus bikes can be relocated to regions with a greater density of commuters, thereby maximizing their utilization.}

Passengers familiar with SM demonstrate a 7\% higher modal share in SM-related modes compared to those unfamiliar, as depicted in Figures \ref{fig:Heuri1000end55000modal_sharesnot familiarwith SM} and \ref{fig:Heuri1000end55000modal_sharesfamiliarwith SM.}.
Those utilizing PT more than once per week exhibit a 5\% higher modal share on combining PT with other modes compared to less frequent users, as evidenced in Figures \ref{fig:Heuri1000end55000modal_sharesuse PT  small onceweek.} and \ref{fig:Heuri1000end55000modal_sharesuse PTlargeonceweek.}. {This indicates that increasing promotional efforts about the availability of SM modes and the benefits of multi-modal transport can potentially boost their usage. Campaigns, information sessions, and enhanced visibility of SM options could be employed to increase familiarity and comfort among potential users, focusing on regions less familiar with SM.}

High-income passengers display a 9\% higher usage of scooters and a 13\% lower usage of PT compared to low-income passengers, as depicted in Figures \ref{fig:Heuri1000end55000modal_shareslowincome.} and \ref{fig:Heuri1000end55000modal_shareshighincome.}.
Female passengers exhibit a 6\% higher usage of scooters and private vehicles, a 6\% lower usage of bikes, as well as a 6\% decrease in combined PT and scooter usage compared to male passengers, as shown in Figures \ref{fig:Heuri1000end55000modal_sharesmale.} and \ref{fig:Heuri1000end55000modal_sharesfemale.}.
Passengers with higher education levels demonstrate a considerably - 22\% - higher share of scooter usage compared to those with lower education levels, as depicted in Figures \ref{fig:Heuri1000end55000modal_sharesloweducation.} and \ref{fig:Heuri1000end55000modal_shareshigheducation.}. {In cases where there is an excess of scooters, a viable strategy would be to move these scooters to areas where there is a higher concentration of high-income and high-education passengers.}

Users of SM exhibit an 8\% higher modal share for scooters and a 7\% lower inclination towards bikes, as shown in Figures \ref{fig:Heuri1000end55000modal_shareshaving used SM.} and \ref{fig:Heuri1000end55000modal_shareshaving notused SM.}. {Businesses and city planners might consider increasing investments in scooter infrastructure where SM usage is high.} Regarding age demographics, passengers under 35 and between 35 and 65 demonstrate similar modal shares. However, passengers over 65 exhibit a substantial - 37\% - lower rate of scooter usage and a 42\% increase in combined PT and other mode usage, alongside a 6\% increase in walking, as illustrated in Figures \ref{fig:Heuri1000end55000modal_sharesage small 35.} to \ref{fig:Heuri1000end55000modal_sharesagelarge65.}. {This shift could be due to mobility challenges, safety concerns, or a preference for less physically demanding modes of transport. Transport services could be tailored to meet the specific needs of different age groups. For the elderly, ensuring safer, more accessible public transport options and facilitating easy transfers between modes could encourage greater independence and mobility.}

\subsubsection{{Analysis of factors influencing multimodal trips}}
\label{Analysis of factors influencing multimodal trips}

\begin{figure}[h]
\centering

\begin{subfigure}[b]{0.395\textwidth}
\includegraphics[width=\textwidth]{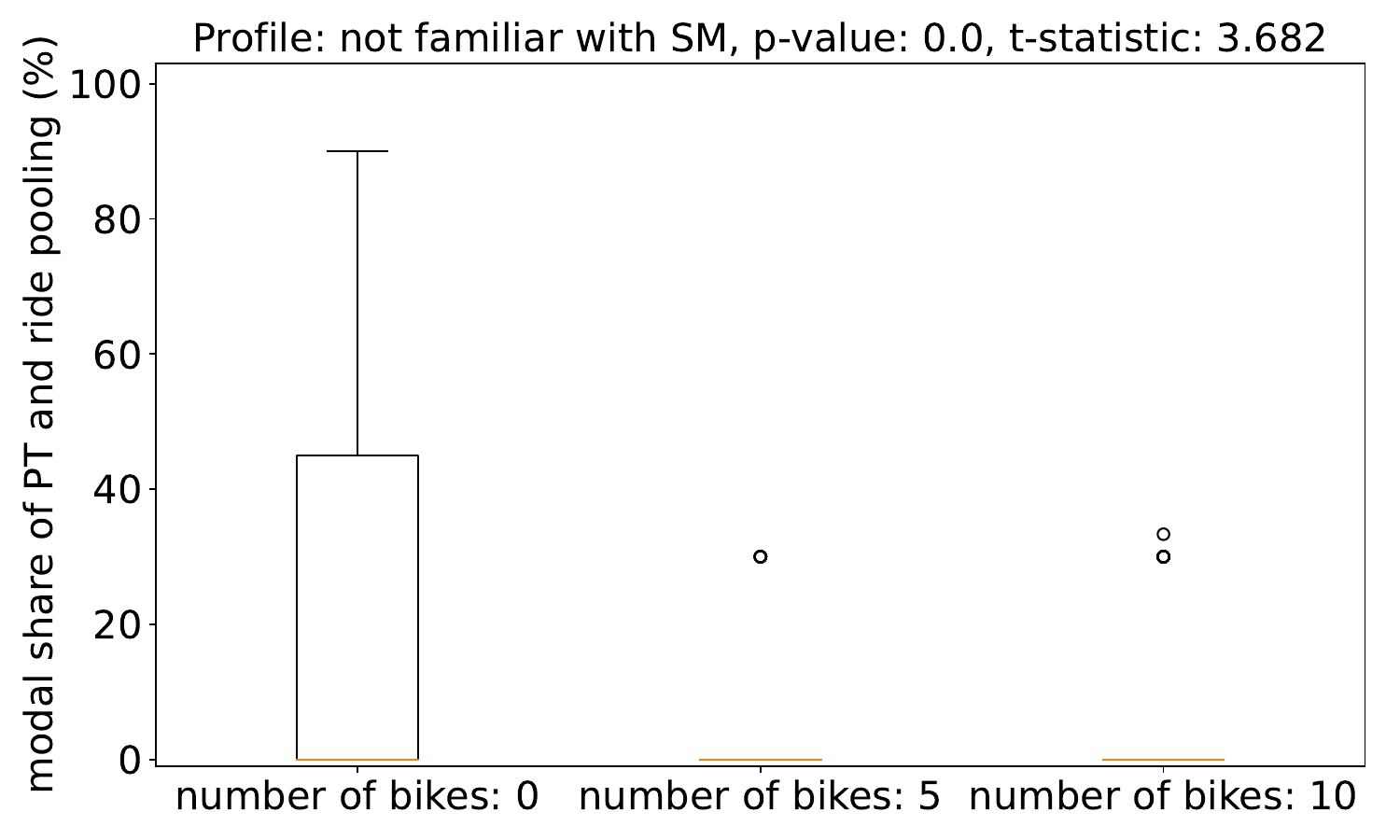}

\caption{}
\label{fig:Heuri1000end55000xnumber of bikesprofile_not familiar with SMy_modal_share_PT_ride}
\end{subfigure}
\begin{subfigure}[b]{0.395\textwidth}
\includegraphics[width=\textwidth]{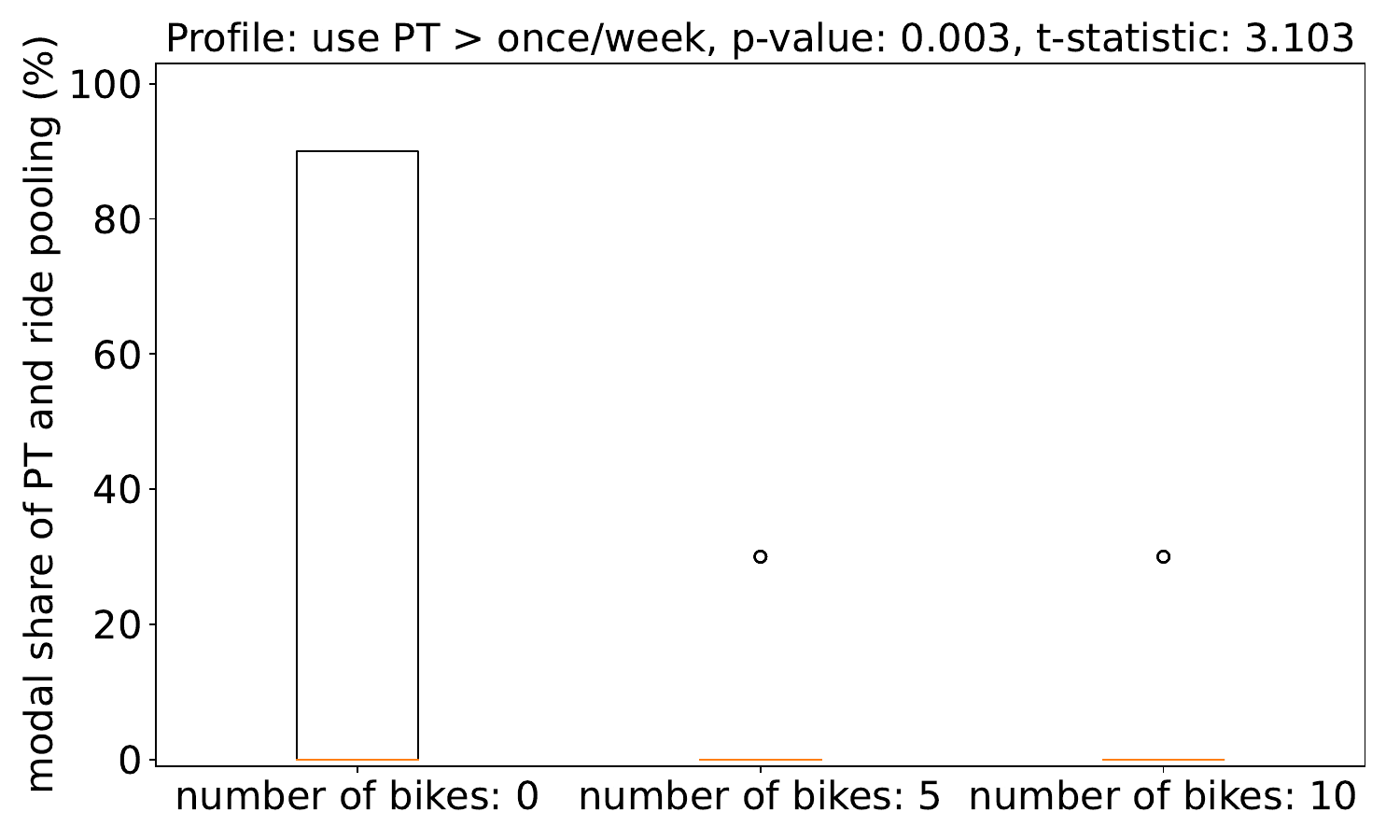}

\caption{}
\label{fig:Heuri1000end55000xnumber of bikesprofile_use PTlargeonceweeky_modal_share_PT_ride}
\end{subfigure}

\begin{subfigure}[b]{0.55\textwidth}
\includegraphics[width=\textwidth]{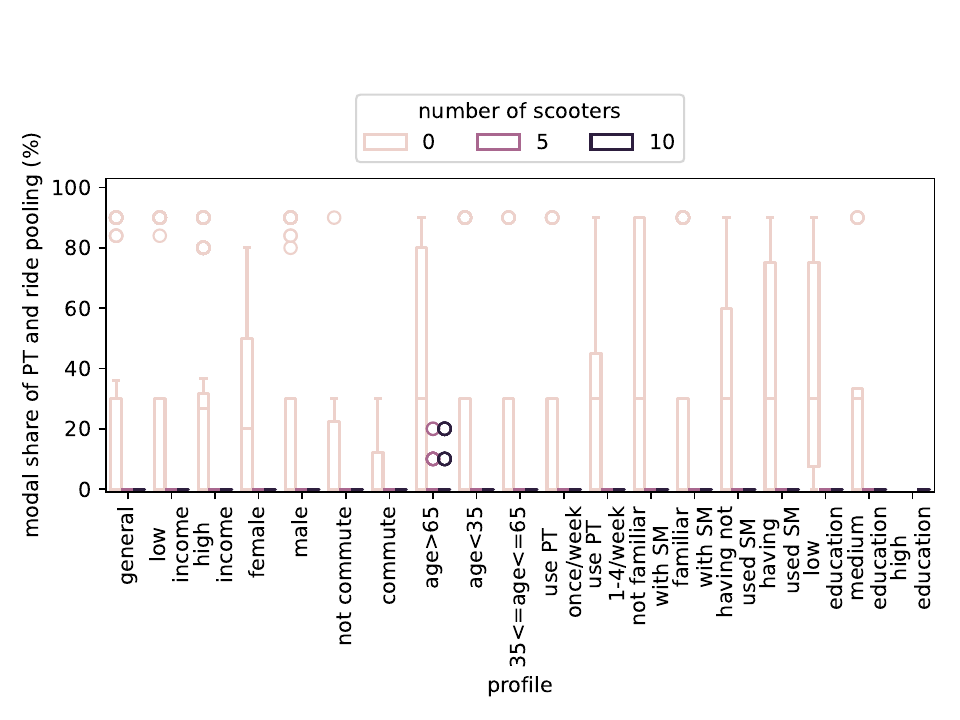}

\caption{}
\label{fig:Heuri1000end55000profile_as_xlegendnumber of scootersx_profiley_modal share of combination of public transport and ride pooling}
\end{subfigure}

\caption{Box plots of the modal share of the combination of PT and ride-pooling}
\label{fig:Modal shares of combinations of PT and ride pooling}
\end{figure}

The modal share of PT combined with ride-pooling is influenced by the availability of micromobility vehicles. For passengers who are not familiar with shared mobility and passengers who use PT frequently, they use PT and ride-pooling when shared bikes are not available. However, the presence of available shared bikes reduces the modal shares of PT and ride-pooling to nearly 0\%, as shown in Figures \ref{fig:Heuri1000end55000xnumber of bikesprofile_not familiar with SMy_modal_share_PT_ride} and \ref{fig:Heuri1000end55000xnumber of bikesprofile_use PTlargeonceweeky_modal_share_PT_ride}. Therefore, the availability and affordability of shared micro-mobility play a crucial role in shaping passengers' preferences for multi-modal transport options. Compared to bikes, the introduction of scooters has a substantial impact on all passenger types, leading to a significant reduction in the modal share of PT combined with ride-pooling, as illustrated in Figure \ref{fig:Heuri1000end55000profile_as_xlegendnumber of scootersx_profiley_modal share of combination of public transport and ride pooling}. {Therefore, increasing the number of bikes and scooters could significantly reduce the use of PT and ride-pooling due to their affordability and convenience for shorter trips. Policymakers and shared micro-mobility providers should increase the availability of shared micromobility vehicles near PT stations to enhance accessibility and encourage integrated transport use.}

\begin{figure}[h]
\centering
\begin{subfigure}[b]{0.395\textwidth}
\includegraphics[width=\textwidth]{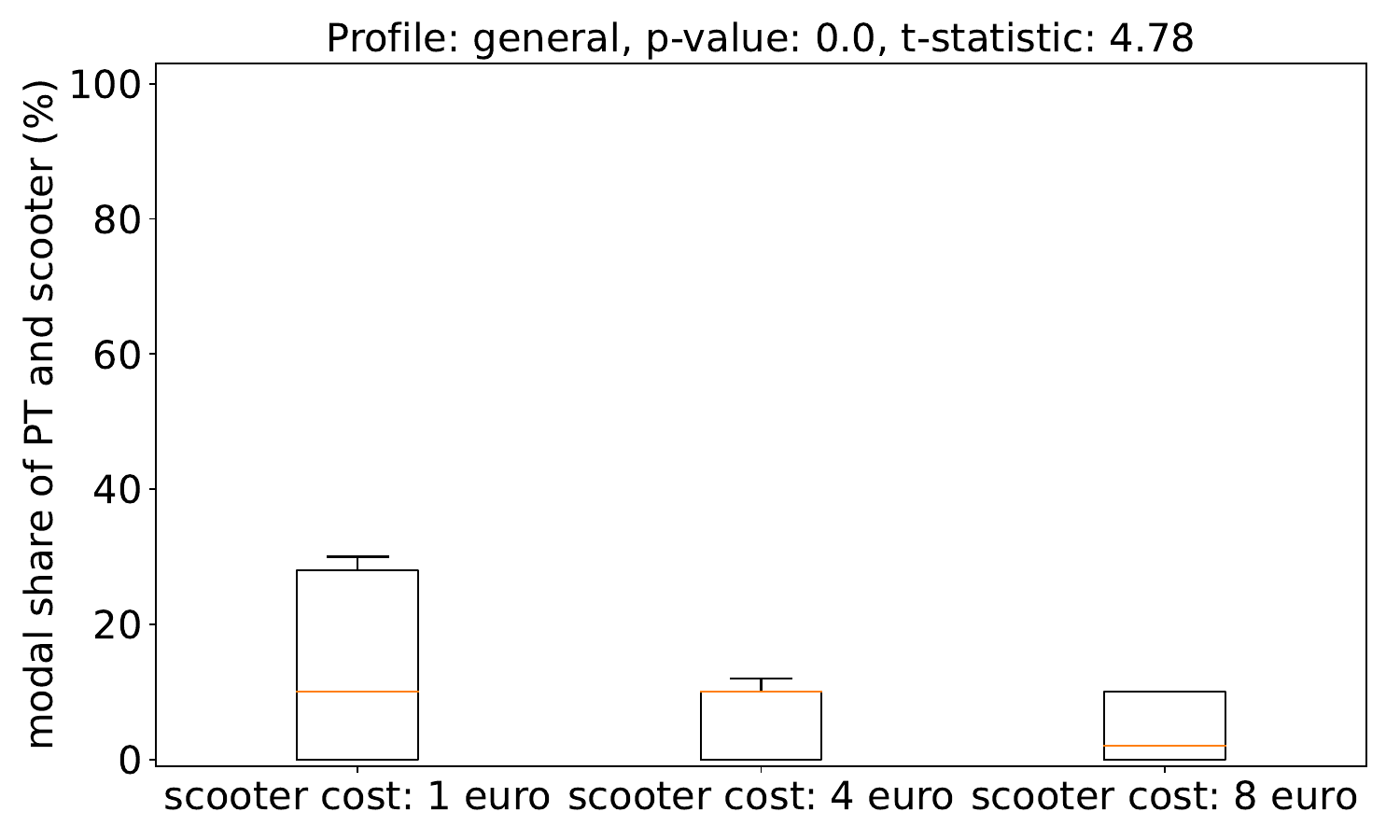}

\caption{}
\label{fig:Heuri1000end55000xscooter costprofile_generaly_modal_share_PT_scooter}
\end{subfigure}
\begin{subfigure}[b]{0.395\textwidth}
\includegraphics[width=\textwidth]{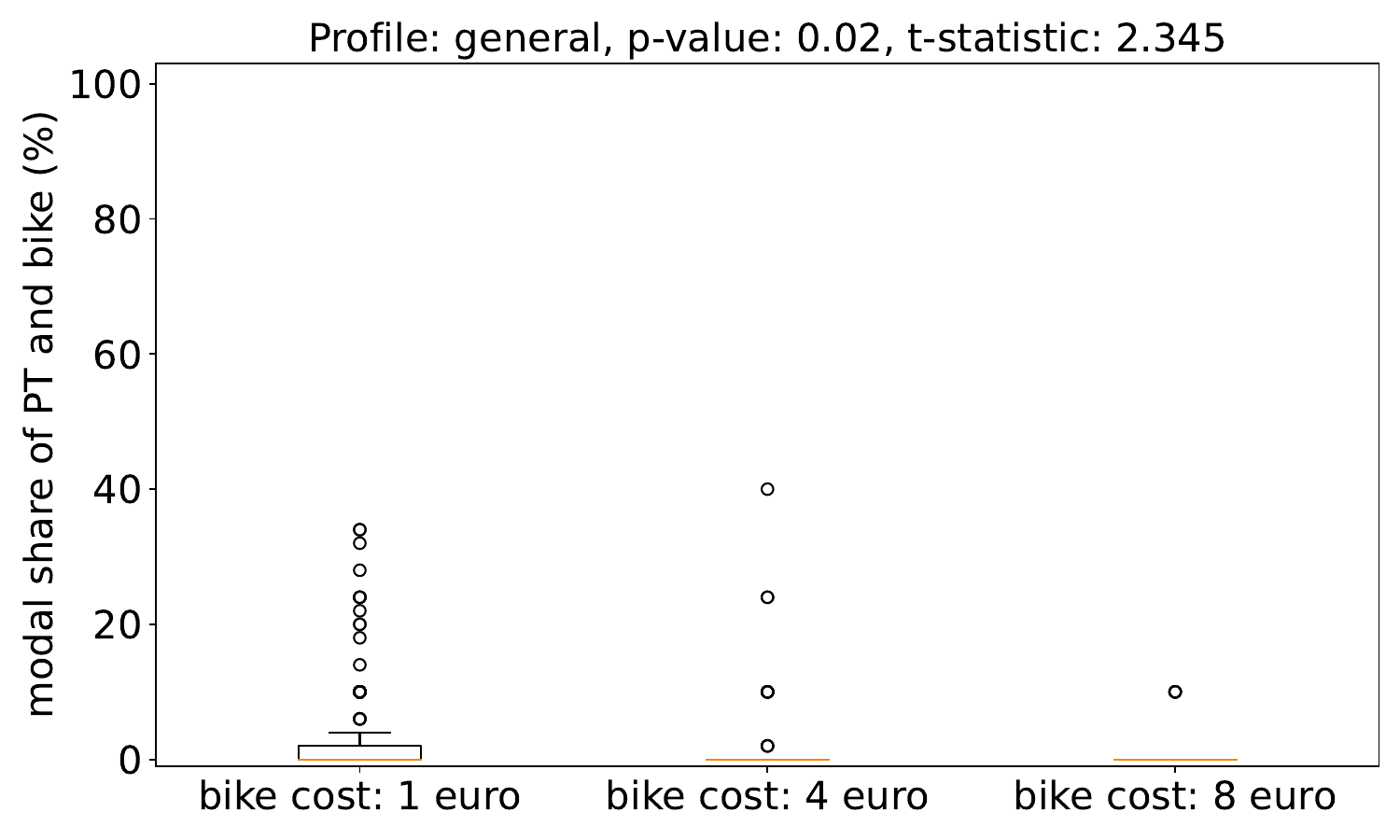}

\caption{}
\label{fig:Heuri1000end55000xbike costprofile_generaly_modal_share_PT_bike}
\end{subfigure}

\begin{subfigure}[b]{0.55\textwidth}
\includegraphics[width=\textwidth]{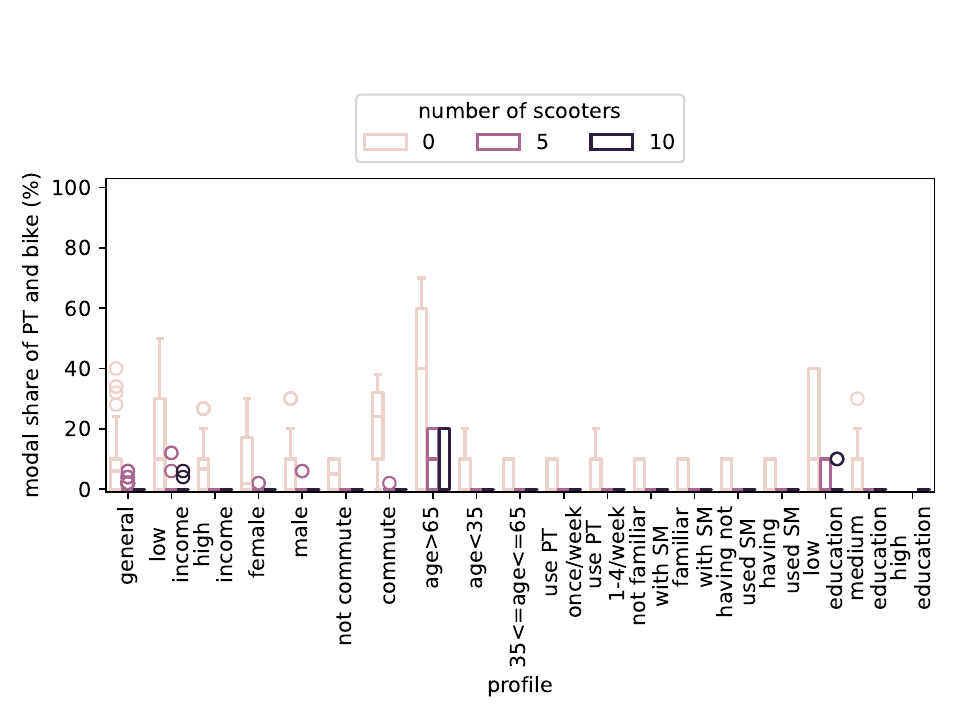}

\caption{}
\label{fig:Heuri1000end55000profile_as_xlegendnumber of scootersx_profiley_modal share of combination of public transport and bike}
\end{subfigure}

\caption{Box plots of the modal shares of a combination of PT and micro-mobility}
\label{fig:Modal shares of combinations of PT and micromobility}
\end{figure}

Figure \ref{fig:Modal shares of combinations of PT and micromobility} depicts the modal shares associated with various combinations of PT and shared mobility. An increase in the costs of scooters and bikes leads to a significant decrease in the modal share of PT and these micromobility options, as shown in Figures \ref{fig:Heuri1000end55000xscooter costprofile_generaly_modal_share_PT_scooter} and \ref{fig:Heuri1000end55000xbike costprofile_generaly_modal_share_PT_bike}. To enhance the modal shares of multi-modal trips, it is essential to keep the cost of shared mobility, particularly for scooters, at a minimum. Affordable scooters will encourage passengers to opt for the combination of PT and scooters, attaining a modal share of 10\%, as illustrated in Figure \ref{fig:Heuri1000end55000xscooter costprofile_generaly_modal_share_PT_scooter}. However, as scooter costs rise, this share drops to 3\%. Figure \ref{fig:Heuri1000end55000profile_as_xlegendnumber of scootersx_profiley_modal share of combination of public transport and bike} shows that the modal share of PT and bikes will be reduced significantly when scooters are available, as scooters substitute bikes. {For businesses in the shared micro-mobility sector, this insight suggests a strategic shift in investment from shared bicycles to scooters, coupled with the necessity to innovate and develop advanced technologies for enhanced management of shared scooter systems.}

\subsubsection{{Analysis of factors influencing single-mode trips}}
\label{Analysis of factors influencing single-mode trips}

\begin{figure}[h]
\centering
\begin{subfigure}[b]{0.395\textwidth}
\includegraphics[width=\textwidth]{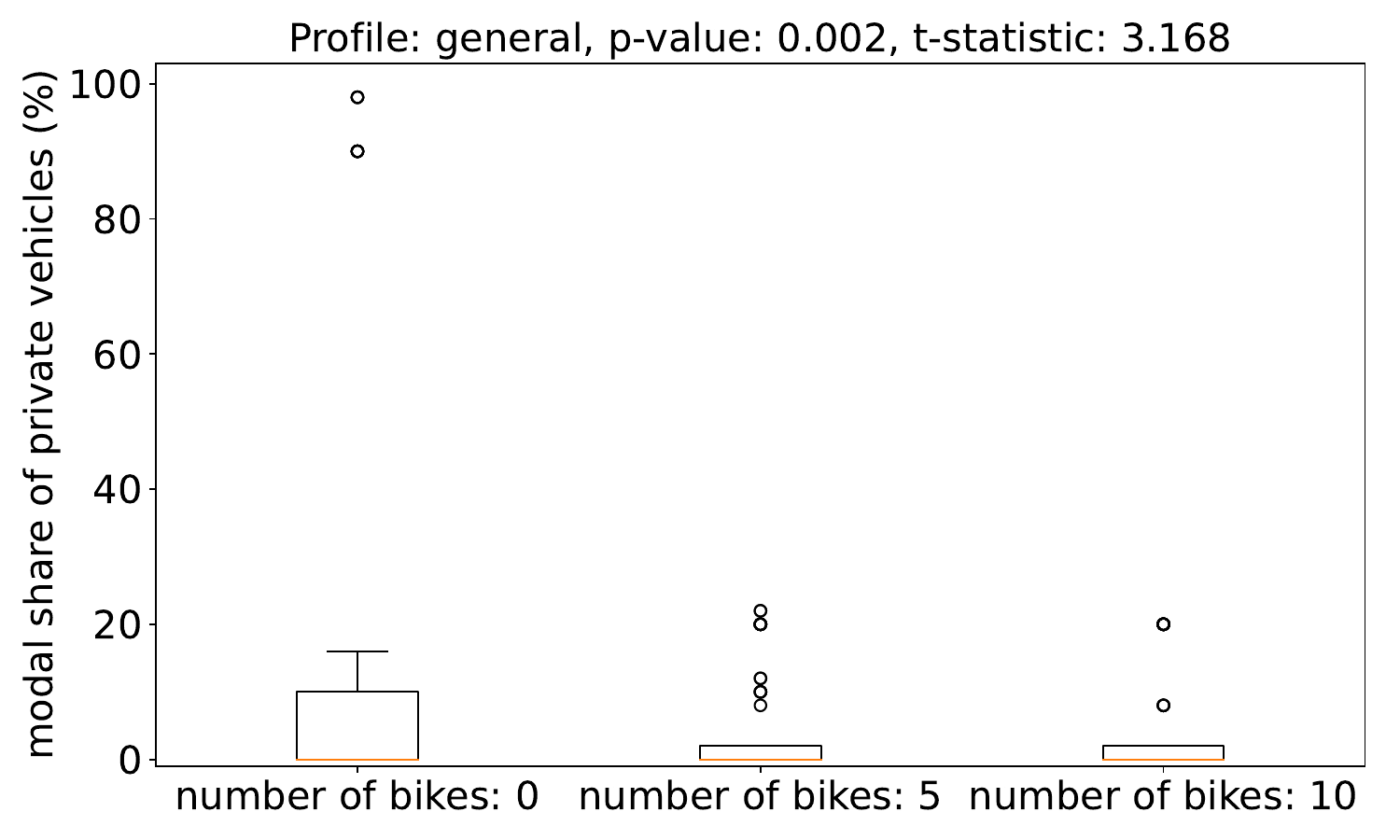}

\caption{}
\label{fig:Heuri1000end55000xnumber of bikesprofile_generaly_modal_share_private}
\end{subfigure}
\begin{subfigure}[b]{0.395\textwidth}
\includegraphics[width=\textwidth]{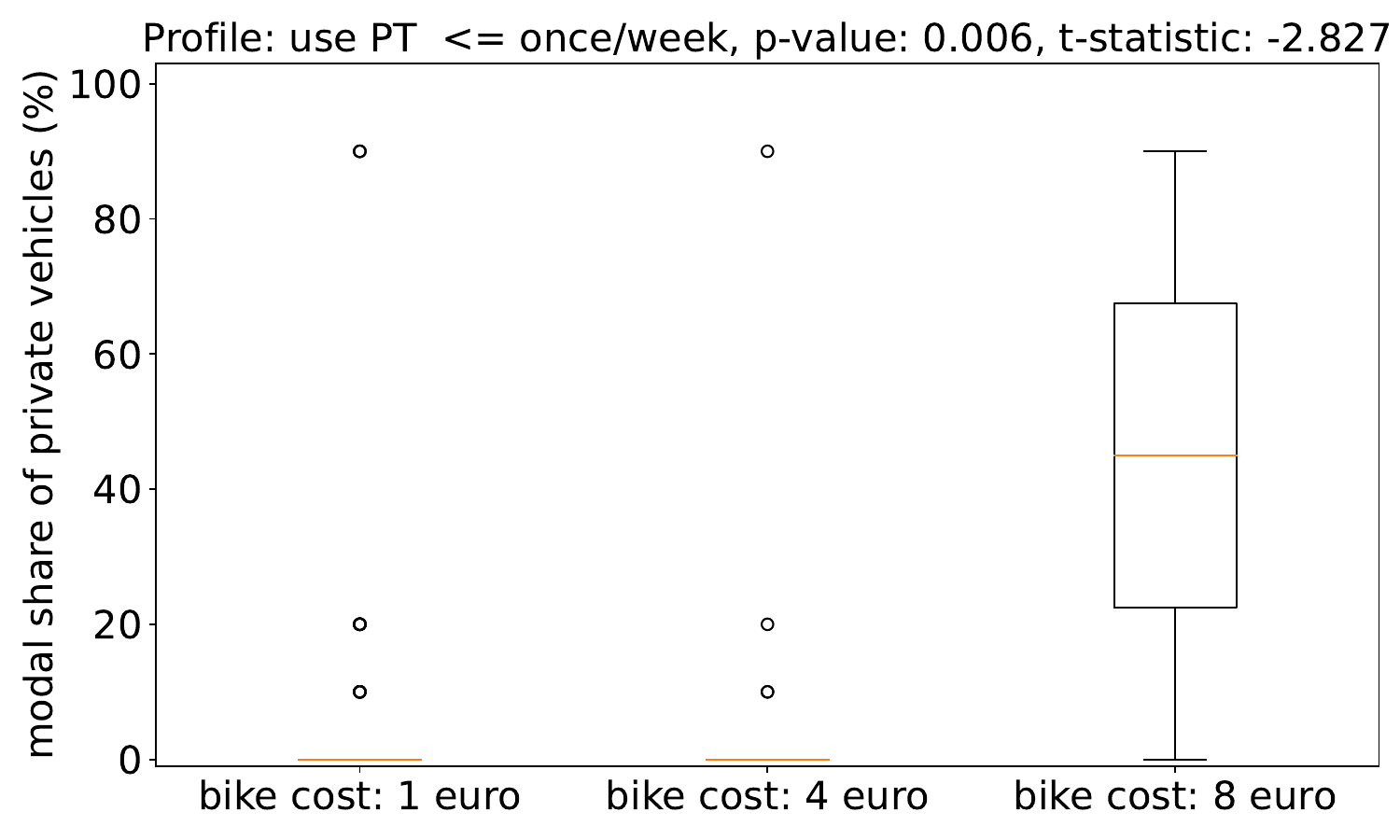}

\caption{}
\label{fig:Heuri1000end55000xbike costprofile_use PT small onceweeky_modal_share_private}
\end{subfigure}
\begin{subfigure}[b]{0.55\textwidth}
\includegraphics[width=\textwidth]{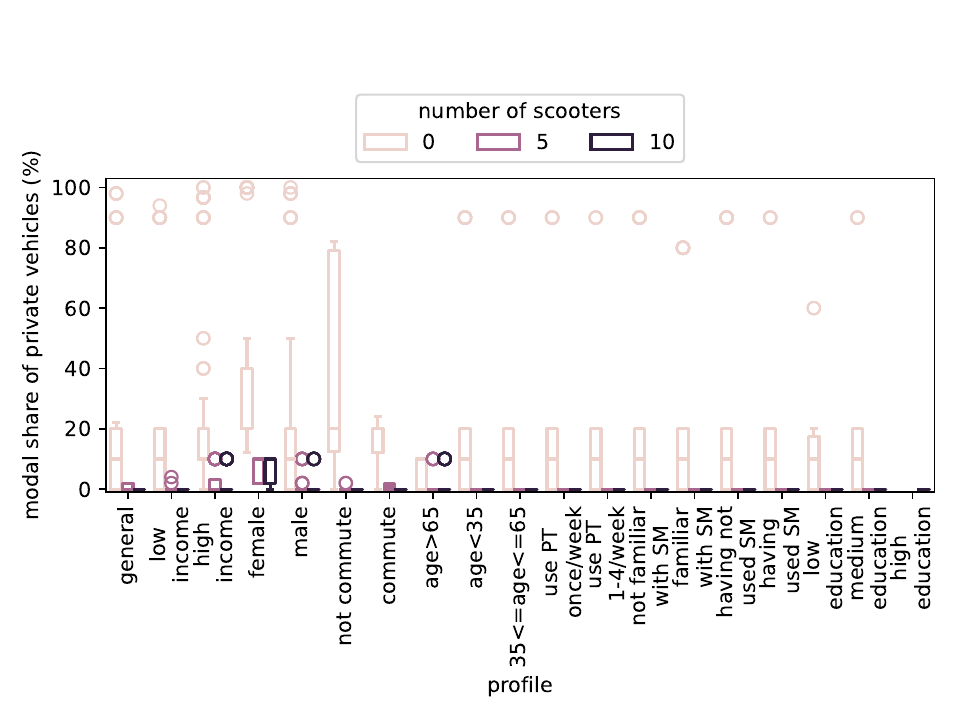}

\caption{}
\label{fig:Heuri1000end55000profile_as_xlegendnumber of scootersx_profiley_modal share of private vehicles}
\end{subfigure}
\caption{Modal shares of private vehicles}
\label{fig:Modal shares of private vehicles}
\end{figure}

Figure \ref{fig:Modal shares of private vehicles} shows the modal share of private vehicles. Increasing the number of available bikes reduces the modal share of private vehicles, as shown in Figure \ref{fig:Heuri1000end55000xnumber of bikesprofile_generaly_modal_share_private}. For most types of passengers, there is no significant change in the modal share of private vehicles when changing the cost of bikes, as these passengers opt for alternative transport modes in response to changes in the cost of renting a shared bike. However, for passengers who use PT once a week or less, when the cost of renting a shared bike becomes too high, there is a significant surge in the use of private vehicles, pushing modal share from 0\% to approximately 50\%, as shown in Figure \ref{fig:Heuri1000end55000xbike costprofile_use PT small onceweeky_modal_share_private}. {Therefore, it is crucial to keep the cost of shared mobility services affordable, particularly for passengers who do not frequently use PT. Finding the right balance between affordability and convenience is key to promoting the widespread adoption of shared mobility options and reducing private vehicle usage effectively.} For all passenger types, the availability of scooters leads to a decreased tendency to use private vehicles, as shown in Figure \ref{fig:Heuri1000end55000profile_as_xlegendnumber of scootersx_profiley_modal share of private vehicles}. {The results indicate that when sufficient affordable services are available, seamless connectivity between SM and PT can be ensured, enhancing passenger convenience and overall social welfare. However, if SM services are limited, resource constraints may force some passengers to rely solely on private vehicles. In such cases, a conflict may arise between individual optimal choices and social welfare, as previously allocated services may leave current passengers without access to SM options.} {To reduce the reliance on private vehicles, the availability of bikes and scooters is imperative; their accessibility typically leads to a decrease in private vehicle usage.} These findings suggest that interventions focused on augmenting micro-mobility availability could serve as effective measures in encouraging sustainable transport choices.

\begin{figure}[h]
\centering
\begin{subfigure}[b]{0.395\textwidth}
\includegraphics[width=\textwidth]{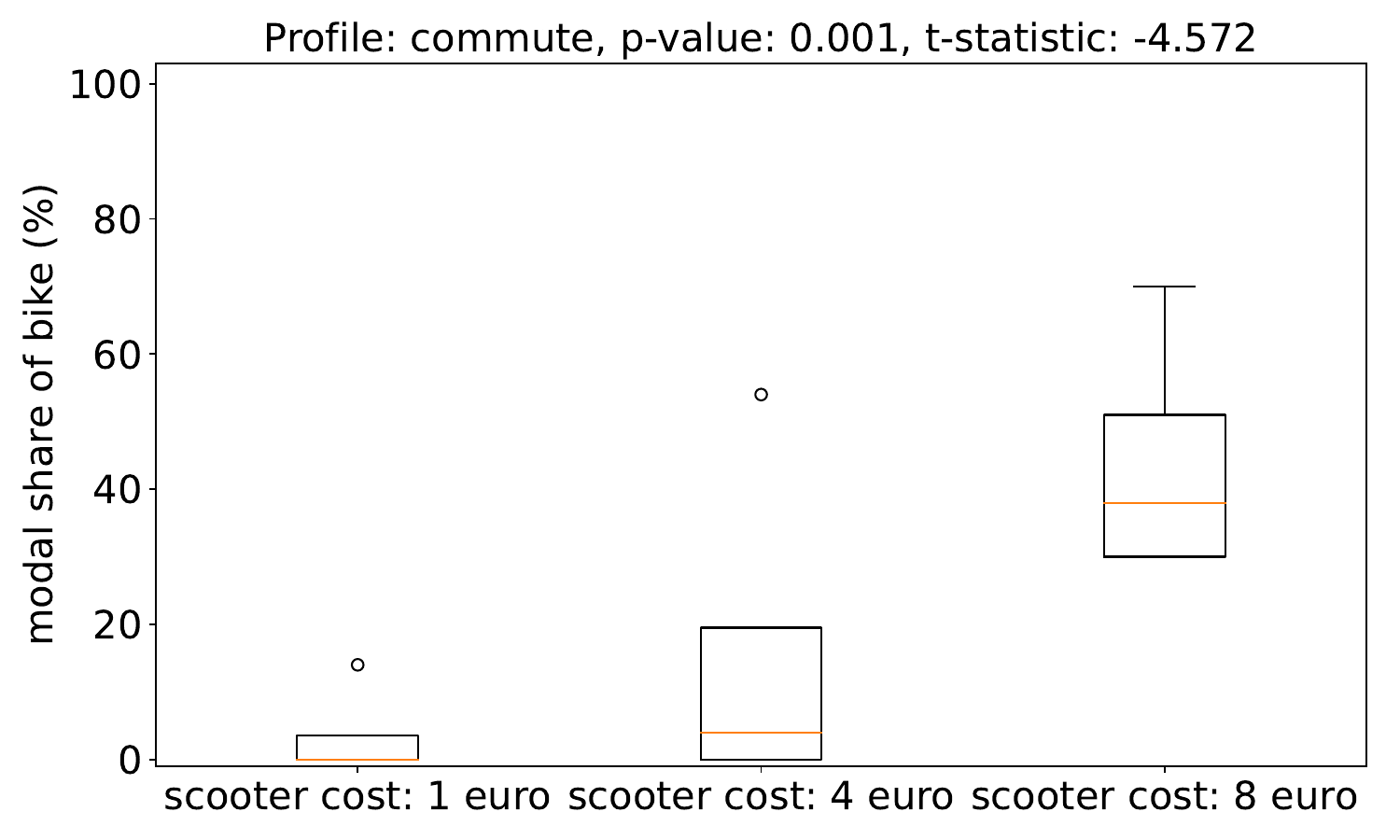}

\caption{}
\label{fig:Heuri1000end55000xscooter costprofile_commutey_modal_share_fixed}
\end{subfigure}
\begin{subfigure}[b]{0.395\textwidth}
\includegraphics[width=\textwidth]{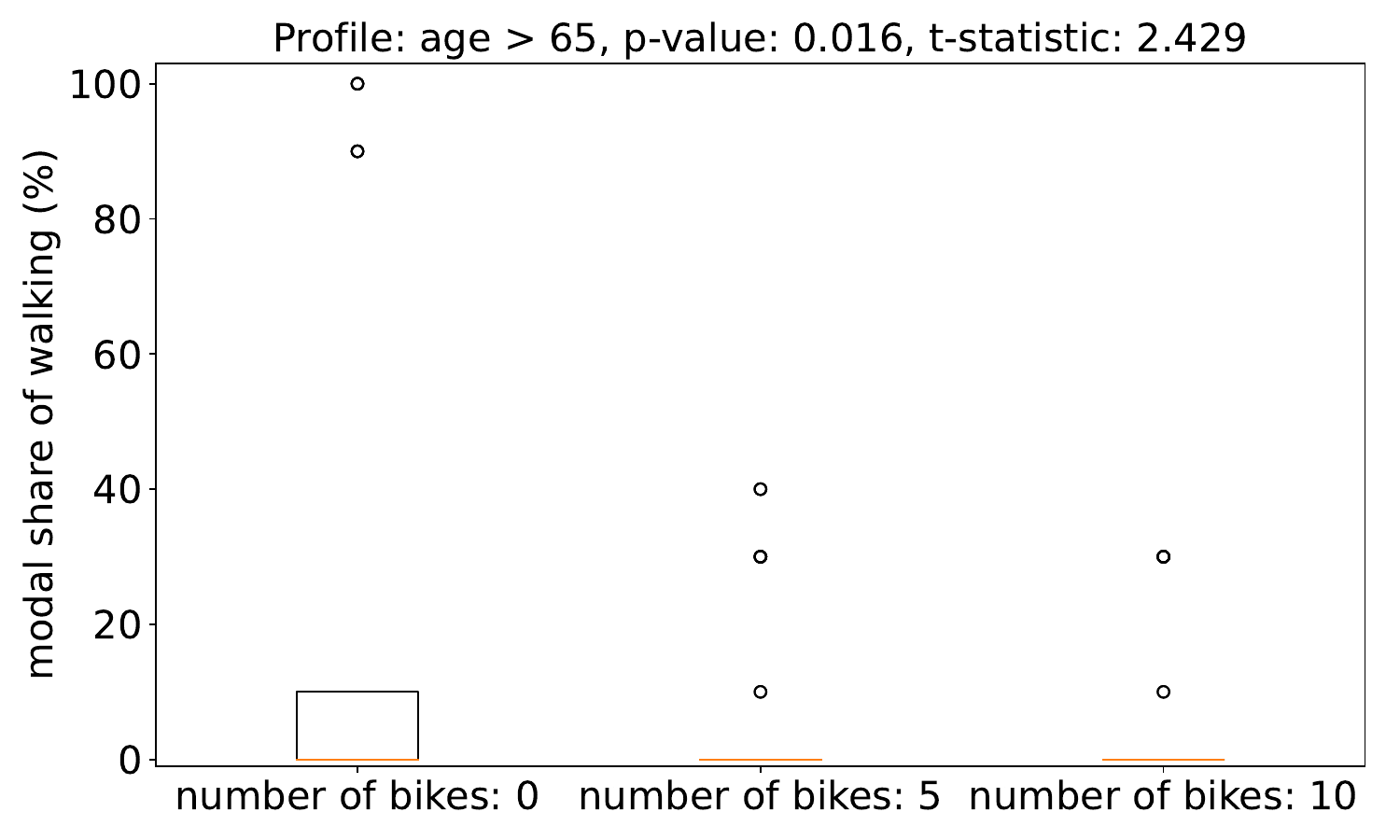}

\caption{}
\label{fig:Heuri1000end55000xnumber of bikesprofile_agelarge65y_modal_share_walking}
\end{subfigure}

\caption{Box plots of the modal shares of bikes and walking}
\label{fig:Modal shares of bikes and walking}
\end{figure}

The costs of shared bikes and scooters influence the modal shares of each other. When the cost of renting a shared scooter increases, the modal share of bikes increases for some user segments, such as commuters, passengers who have used shared mobility, and passengers with low and medium education levels. Figure \ref{fig:Heuri1000end55000xscooter costprofile_commutey_modal_share_fixed} illustrates the case of commuters, and the t-statistic is -4.572 when contrasting scooter costs of 1 euro and 8 euros. {The increase in scooter prices impacts commuters, resulting in a notable shift of 40\% of the modal share from scooters to bikes.} For individuals aged over 65, increasing the number of shared bikes can significantly reduce their walking time, suggesting that shared bikes offer a time-saving alternative for old people, as shown in Figure \ref{fig:Heuri1000end55000xnumber of bikesprofile_agelarge65y_modal_share_walking}.

\subsubsection{{Analysis of factors influencing travel time and waiting time}}
\label{Analysis of factors influencing travel time and waiting time}

\begin{figure}[h]
\centering
\begin{subfigure}[b]{0.295\textwidth}
\includegraphics[width=\textwidth]{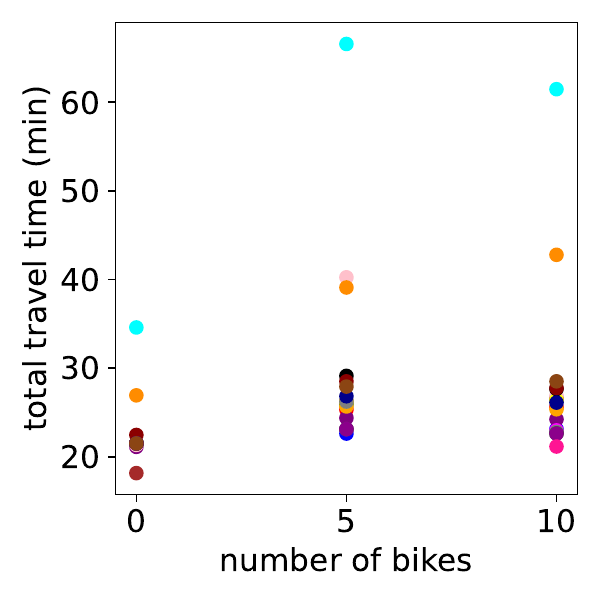}
\caption{}
\label{fig:Heuri1000end55000travel_timestotal travel timenumber of bikes}
\end{subfigure}
\begin{subfigure}[b]{0.295\textwidth}
\includegraphics[width=\textwidth]{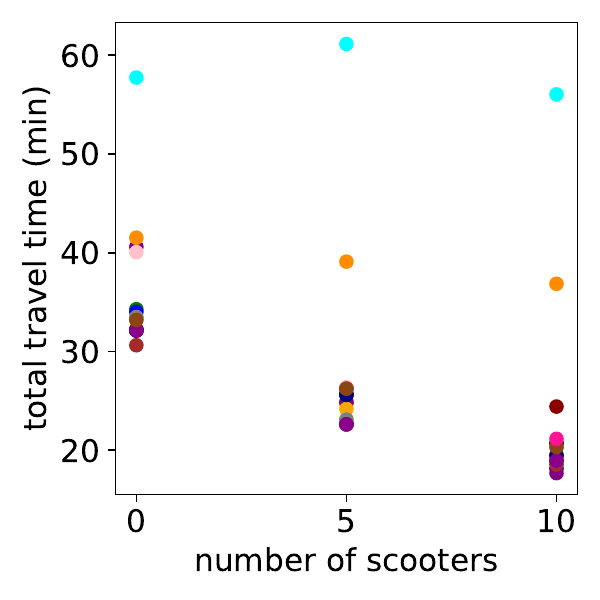}
\caption{}
\label{fig:Heuri1000end55000travel_timestotal travel timenumber of scooters}
\end{subfigure}
\begin{subfigure}[b]{0.295\textwidth}
\includegraphics[width=\textwidth]{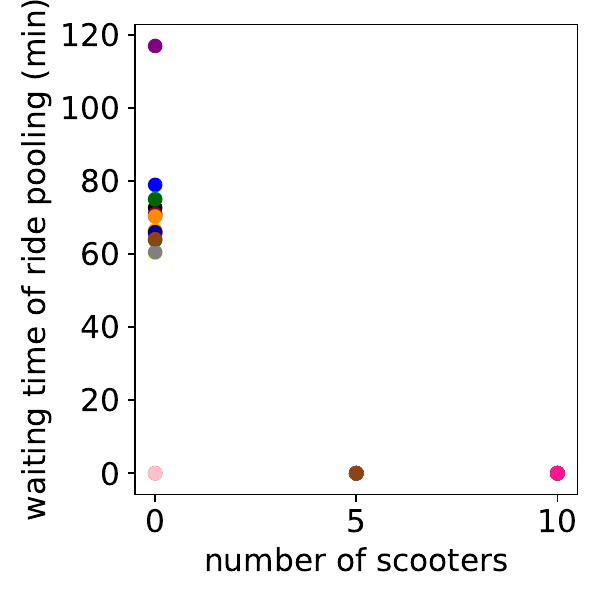}
\caption{}
\label{fig:Heuri1000end55000travel_timeswaiting time of ride pooling (min)number of scooters}
\end{subfigure}

\begin{subfigure}[b]{0.295\textwidth}
\includegraphics[width=\textwidth]{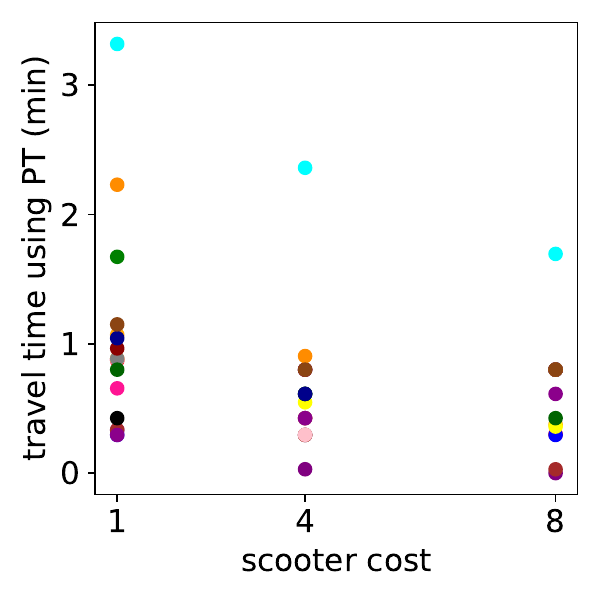}
\caption{}
\label{fig:Heuri1000end55000travel_timestravel time using PT (min)scooter cost}
\end{subfigure}
\begin{subfigure}[b]{0.295\textwidth}
\includegraphics[width=\textwidth]{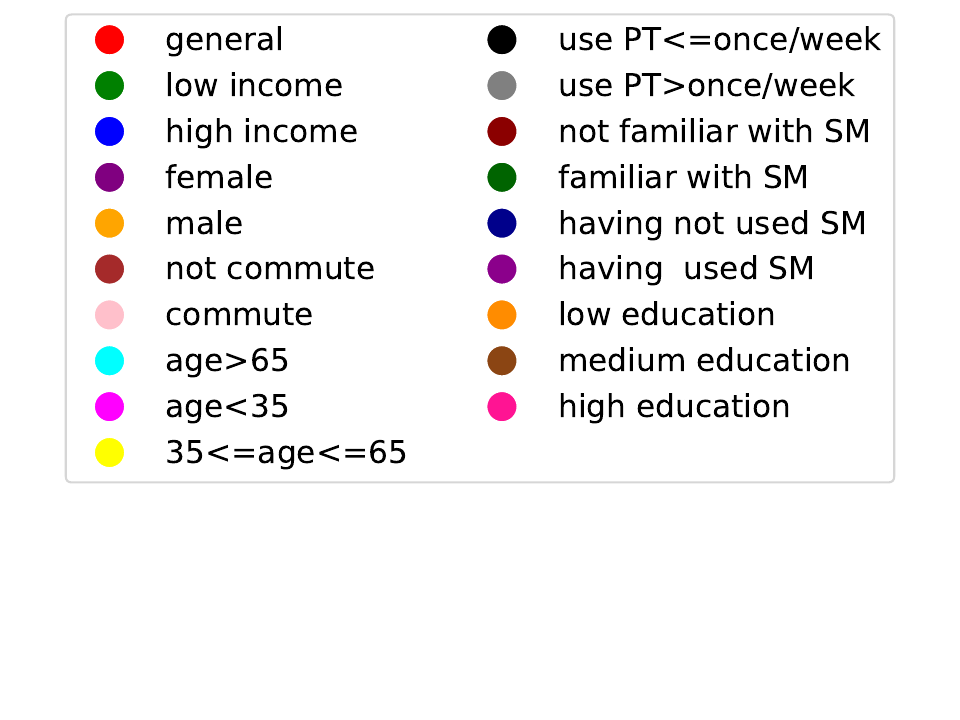}
\label{fig:Heuri1000end55000travel_times_only_legends}
\end{subfigure}
\caption{Travel time and waiting time}
\label{fig:Travel time, idle time, and waiting time}
\end{figure}

Figure \ref{fig:Travel time, idle time, and waiting time} illustrates the variations in travel time and waiting time across different scenarios. The impact of increasing the size of the shared bike fleet on total travel time is not as pronounced as the effect of increasing the size of the shared scooter fleet, as shown in Figure \ref{fig:Heuri1000end55000travel_timestotal travel timenumber of bikes} and \ref{fig:Heuri1000end55000travel_timestotal travel timenumber of scooters}. The total travel time for almost all passenger types decreases with an increase in the size of the scooter fleet due to scooters' speed advantage. {Scooters, especially when integrated with public transport, offer a faster alternative compared to bikes.} Furthermore, encouraging the deployment of scooters can contribute to overall time savings for various passenger demographics. Increasing the cost of scooters reduces the waiting time for ride pooling, as shown in Figure \ref{fig:Heuri1000end55000travel_timeswaiting time of ride pooling (min)number of scooters}. This is because fewer passengers opt for ride pooling when scooters are readily available. Interestingly, Figure \ref{fig:Heuri1000end55000travel_timestravel time using PT (min)scooter cost} demonstrates that scooter cost even influences the travel time when using PT. When scooters are affordable, more passengers opt for the combination of PT and scooters, but as scooter costs rise, passengers use less PT. Therefore, reducing scooter costs to encourage combined PT and scooter usage may contribute to a more integrated and efficient transport system.

\section{Conclusion}\label{Conclusions}

{To reduce emissions by promoting public transport and shared mobility services as alternatives to private vehicles,} we developed a comprehensive mathematical model alongside a tailored heuristic algorithm designed specifically for dynamic preference-based multi-modal trip planning. This approach takes into account the diverse preferences of passengers within a dynamic operational environment, incorporating multi-modal travel constraints, such as time windows, capacity of vehicles, and transfers. We explicitly account for different user segments in terms of socio-demographic and travel characteristics. Multi-modal planning includes the consideration of line- and schedule-based public transport services as well as shared mobility in the form of ride-pooling services, shared bikes, and shared scooters. We compared the results of the heuristic algorithm with those obtained from the exact approach. Our findings reveal that the proposed meta-heuristic algorithm demonstrates remarkable efficiency in achieving solutions close to optimality within significantly reduced computation times.

{The developed meta-heuristic algorithms scale these solutions to a real case in a suburban area of Rotterdam, which not only validates our approach but also demonstrates its potential value in addressing real-world transportation challenges. Considering passengers' preferences is pivotal for informing transport planning and policy decisions for fostering sustainable and efficient transport modes.} For instance, the notable decrease in scooter usage among individuals over 65, coupled with an increase in combined PT and other mode usage, highlights the need for targeted initiatives. Tailored efforts, like promoting public transport for non-commuters or mitigating obstacles to scooter usage among older passengers, can be developed leveraging these insights. Strategies to encourage the adoption of PT services, especially among specific demographics like high-income and female passengers, should also be considered.

{In summary, this study offers actionable insights for transport managers and policymakers, providing a robust framework that supports more informed decision-making and the development of efficient, adaptable urban mobility systems.} {While the current study validates the proposed framework using simulated scenarios, the future study will involve testing the framework in a larger transport network, further evaluating its scalability and applicability in real-world settings.} {The framework allows for the integration of additional shared mobility solutions, such as carsharing, or the exclusion of specific transport modes modeled in this study, such as ride-pooling services. Key parameters, including vehicle capacities, cost structures, and user preferences, can be adjusted to align with the unique characteristics of new mobility modes and adapt to the distinct needs of suburban or urban environments.}

The following research directions are suggested for future research. First, in multi-modal trip planning, the MaaS platform can incentivize passengers to relocate vehicles to areas where they are in higher demand; thus, the service can better meet the needs of users. Therefore, developing passenger-centric re-balancing algorithms is a promising research direction. Machine learning techniques can be utilized to analyze historical trip data and predict demand patterns for different locations and time periods, informing the re-balancing algorithms to devise more efficient strategies. Second, pricing models can be designed to tailor to the preferences of different passenger segments. Dynamic pricing strategies that adjust fares based on factors such as time of day, route popularity, and vehicle availability can be developed, ensuring optimal revenue generation while still meeting passenger demand. Third, future research could employ the Mixed Logit model to better account for taste heterogeneity and correlated alternatives, though this would introduce non-linearity in optimization and require advanced methods to handle the increased complexity. {Fourth, when serving a batch of requests at each time step, conflicts may arise between user-optimal and system-optimal solutions. In such cases, passengers can be nudged to achieve an overall more efficient transport. Therefore, future research can develop a variety of nudging strategies tailored to different passenger segments based on their preferences. A menu of multi-modal alternatives that incorporates these nudges can be designed, including options such as discounted fares during off-peak hours, loyalty rewards for frequent users, or personalized recommendations based on past travel behavior. Last, future research could incorporate dynamic preferences that evolve with network conditions, offering a more accurate representation of real-world passenger behavior over time.}

\ACKNOWLEDGMENT{This research is supported by the project ``Seamless Shared Urban Mobility (SUM)". This project has received funding from the European Union's Horizon 2020 research and innovation program under Grant Agreement no. 101103646.}

\bibliographystyle{dcu}
\bibliography{manuscript}

\renewcommand{\thetable}{A.\arabic{table}}
\setcounter{table}{0}
\begin{APPENDICES}

\section{Parameters in the utility function}
\label{Parameters in the utility function}

Tables \ref{table:utility parameters} and \ref{table:utility parameters_beta} show the Alternative specific constant (ASC) and $\beta$ parameters for different profiles in the utility function. {The survey was conducted using a stated preference (SP) approach, employing an online questionnaire distributed to participants in the Rotterdam area and neighboring municipalities. The hypothetical scenarios assumed fully integrated systems of shared micromobility and public transport. Data on attributes like travel time and cost were collected alongside socio-demographic and transport-related characteristics such as age, income, frequency of public transport use, previous familiarity with shared micromobility, etc. These factors are detailed under the ``Profile" column in Tables \ref{table:utility parameters} and \ref{table:utility parameters_beta}. The preference parameters for passengers with different profiles are presented in the columns labeled ``ASC" and ``$\beta$" in Tables \ref{table:utility parameters} and \ref{table:utility parameters_beta}.}

\begin{table}[!hp]
\centering

\begin{adjustbox}{max width=\textwidth}
\begin{threeparttable}
\caption{ASC parameters for different profiles in the utility function.}
\label{table:utility parameters}
\begin{tabular}{l c c c c c c c c}

 \toprule

  Profile&$ASC_{\text{metro}}$& $ASC_{\text{bus\&tram}}$& $ASC_{\text{sb}}$& $ASC_{\text{ss\&sb\&ride}}^{\text{sub}}$& $ASC_{\text{walk}}$& $ASC_{\text{ss}}$& $ASC_{\text{bike}}$& $ASC_{\text{car\&ride}}$ \\
\midrule
A. aggregated &-0.865&0.683&-0.856&-0.934&0.007&-1.380&-0.275&0\\
  B1. commute&-0.120&0.300&0.437&0.173&0.173&0.036&0.048&0\\
  B2. not commute&-0.796&0.524&-1.090&-1.030&0.173&-1.400&-0.297&0\\
  C1. familiar with SM&-0.851&0.687&0.187&0.439&0.007&0.856&-0.274&0\\
  C2. not familiar with SM&-0.851&0.687&-0.986&-1.300&0.007&-2.110&-0.274&0\\
D1. use PT $\leqslant$ once/week&-1.420&0.697&-0.894&-1.130&0.007&-1.740&-0.479&0\\
D2. use PT $\geq$ once/week&1.270&-0.186&-0.263&0.188&0.007&0.839&0.498&0\\
E1. low income&-0.687&0.664&-0.773&-0.815&0.007&-1.310&-0.357&0\\

 E2. high income&0.163&-0.683&-0.536&-0.525&0.007&0.341&0.238&0\\
  F1. male&-0.496&0.157&-0.341&0.603&-0.210&-0.210&-0.265&0\\
F2. female&-0.642 &0.621 &-0.743 &-1.210 &-0.210 &-1.290 &-0.153 &0\\
G1. low education&-0.836&1.140&-0.773&-1.090&0.007&-2.030&-0.555&0\\
G2. high education&0.342&-0.907&-0.223&0.016&0.007&0.016&0.638&0\\
H2. have not used SM&-0.866&0.684&-0.993&-1.200&0.007&-1.690&-0.275&0\\
H1. used SM&-0.866&0.684&0.659&1.080&0.007&1.150&-0.275&0\\
 I3. age $\geqslant$ 65&-0.606&1.130&-1.970&-1.490&0.007&-2.320&-0.710&0\\
I1. age $\leqslant$ 35&-0.275&-0.836&1.470&0.836&0.007&1.290&0.661&0\\
I2. 35 $\geq$ age $\leq$ 65& -0.323&-0.541&1.240&0.655&0.007&0.957&0.424&0\\
\bottomrule

\end{tabular}

  \end{threeparttable}
  \end{adjustbox}

\end{table}

\vspace{-0.5cm}

\begin{table}[!hp]
\centering

\begin{adjustbox}{max width=\textwidth}
\begin{threeparttable}
\caption{$\beta$ parameters for different profiles in the utility function.}
\label{table:utility parameters_beta}
\begin{tabular}{l c c c c c c}

 \toprule

  Profile& $\beta_{\text{mainCost}}$& $\beta_{\text{subCost}}$& $\beta_{\text{mainTime}}$& $\beta_{\text{subTime}}$& $\beta_{\text{PTwait}}$& $\beta_{\text{walk}}$\\
\midrule
A. aggregated &-0.093&-0.425&-0.034&-0.039&-0.014&-0.064\\
  B1. commute&0.001&0.104&-0.006&0.000&0.035&-0.023\\
  B2. not commute&-0.094&-0.481&-0.031&-0.039&-0.031&-0.053\\
  C1. familiar with SM&-0.035&-0.246&-0.034&-0.039&-0.014&-0.064\\
  C2. not familiar with SM&-0.066&-0.259&-0.034&-0.039&-0.014&-0.064\\
D1. use PT $\leqslant$ once/week&-0.081&-0.668&-0.019&-0.048&-0.074&-0.064\\
D2. use PT $\geq$ once/week&-0.043&0.268&-0.027&-0.002&0.111&0.003\\
E1. low income&-0.114&-0.279&-0.035&-0.024&-0.003&-0.064\\

 E2. high income&0.004&-0.348&-0.004&-0.063&0.011&-0.064\\
  F1. male&0.032&-0.037&0.011&-0.001&-0.005&-0.006\\
F2. female&-0.108 &-0.411 &-0.038 &-0.039 &-0.013 &-0.061\\
G1. low education&-0.068&-0.189&-0.032&0.016&-0.070&-0.030\\
G2. high education&-0.061&-0.451&-0.011&-0.114&0.102&-0.046\\
H2. have not used SM&-0.066&-0.411&-0.034&-0.039&-0.014&-0.064\\
H1. used SM&-0.035&-0.246&-0.034&-0.039&-0.014&-0.064\\
 I3. age $\geqslant$ 65&-0.065&-0.131&-0.020&0.043&-0.056&-0.009\\
I1. age $\leqslant$ 35&-0.046&-0.502&-0.028&-0.146&0.064&-0.093\\
I2. 35 $\geq$ age $\leq$ 65&-0.031&-0.381&-0.014&-0.102&0.044&-0.064\\
\bottomrule

\end{tabular}

  \end{threeparttable}
  \end{adjustbox}

\end{table}

\end{APPENDICES}

\end{document}